# LARGE DEVIATION ASYMPTOTICS FOR OCCUPANCY PROBLEMS[1]

By Paul Dupuis[2], Carl Nuzman and Phil Whiting

*Brown University, Bell Labs and Bell Labs*

In the standard formulation of the occupancy problem one considers the distribution of $r$ balls in $n$ cells, with each ball assigned independently to a given cell with probability $1/n$. Although closed form expressions can be given for the distribution of various interesting quantities (such as the fraction of cells that contain a given number of balls), these expressions are often of limited practical use. Approximations provide an attractive alternative, and in the present paper we consider a large deviation approximation as $r$ and $n$ tend to infinity. In order to analyze the problem we first consider a dynamical model, where the balls are placed in the cells sequentially and "time" corresponds to the number of balls that have already been thrown. A complete large deviation analysis of this "process level" problem is carried out, and the rate function for the original problem is then obtained via the contraction principle. The variational problem that characterizes this rate function is analyzed, and a fairly complete and explicit solution is obtained. The minimizing trajectories and minimal cost are identified up to two constants, and the constants are characterized as the unique solution to an elementary fixed point problem. These results are then used to solve a number of interesting problems, including an overflow problem and the partial coupon collector's problem.

**1. Introduction.** Urn occupancy problems center on the distribution of $r$ balls in $n$ cells, typically with each ball independently assigned to a given cell with probability $1/n$. The literature on the general topic is enormous. See, for example, [9, 19, 20] and the references therein.

Received April 2002; revised October 2003.
[1]Supported in part by NSF Grant DMS-03-06070.
[2]Supported in part by NSF Grants DMS-00-72004, ECS-99-79250 and Army Research Office Grants DAAD19-00-1-0549 and DAAD19-02-1-0425.

There are many different questions one can pose. For example, it may be that one is interested in the distribution of $(\Gamma_0, \Gamma_1, \dots)$, where $\Gamma_i$ is the number of cells containing exactly $i$ balls after all $r$ balls have been thrown. In the *classical* occupancy problem [9, 14, 20], one is interested only in the distribution of unoccupied urns $\Gamma_0$. In other cases, one might be interested in the (random) number of balls required to fill all cells to a given level, or the number required so that a given fraction are filled to that level—the so called coupon collector's or dixie cup problem. In biology the inverse problem of estimating the number of balls thrown from the number of occupied cells also arises and is used to estimate species abundance [17]. Related applications in biology appear in [2, 3], and an application in computer science appears in [21].

Another related problem of interest is the *overflow problem*, in which the urns are supposed to have a finite capacity $C$, and the number of balls that overflow is the random variable of interest. Ramakrishna and Mukhopadhyay [24] describes an application in computer science concerned with memory access, and [12] considers an application to optical switches. In [12] one is concerned with dimensioning the number of wavelength converters so as to reduce the probability of packet loss across the switch to an acceptable level. In [12] any-color-to-any-color converters are considered. However, by extending the results proved here to the case where the balls have distinct colors, dimensioning in the case where we have many-to-one color converters can be carried out and the blocking probability of the output estimated.

A wide range of results have been proved for the occupancy problem by using "exact" approaches. For example, combinatorial methods are used in [9, 14, 20], and methods that utilize generating functions are discussed in [20]. The implementation of these results, however, can be difficult. For example, in applying the combinatorial results one must compute the difference between large quantities that appear in the inclusion–exclusion formula for the probability that a given fraction of cells are occupied. An analogous difficulty occurs with techniques based on moment generating functions, since one must invert the generating function itself.

Asymptotic methods provide an attractive alternative to both of these approaches. One reason is that they often offer good approximations with only a modest computational effort. A second, perhaps more important reason, is that superior qualitative insights can often be obtained. Indeed, a range of asymptotic results have already been obtained for these problems (see, e.g., [18]). The first large deviations principle (LDP) for urn problems that we are aware of was established in [28] for the special case of the classical occupancy problem. This result was applied in [21] to a boolean satisfiability problem in computer science. Reference [6], which appeared while the present paper was under review, proves a LDP for the infinite-dimensional occupancy measures associated with occupancy processes. The present paper focuses on finite-dimensional occupancy measures in which urn occupancies above a given level are not distinguished. In this finite case, we are able to provide a concise large deviations proof along with explicit, insightful and computable expressions for the rate functions and for the large deviations extremals. The rate function for the occupancy model after all the balls have been thrown is shown to have a simple and fairly explicit rate function, which can be defined in terms of relative entropy with respect to the Poisson distribution. Many different problems can be solved in this framework simply by changing the set over which the rate function is minimized. We also give sample path results for the evolution of the urn occupancies toward a particular event. In principle, the explicit form of the minimizing trajectories for sample path results should enable accurate empirical estimates of unlikely events using importance sampling; see [7, 22].

Let $\lfloor a \rfloor$ denote the integer part of a scalar $a$. With $n$ cells available and a total of $r = \lfloor \beta n \rfloor$ balls to be thrown (with $\beta > 0$), we consider the large deviation asymptotics as $n \to \infty$. The precise statement is as follows. Fix a positive integer $I$. Then with $\Gamma_i^n(\beta)$ equal to the fraction of cells containing $i$ balls, we characterize the large deviation asymptotics of the random vectors $\{(\Gamma_0^n(\beta), \Gamma_1^n(\beta), \ldots, \Gamma_I^n(\beta)), n = 1, 2, \ldots\}$ as $n \to \infty$. A direct analysis of this problem is difficult, and, in fact, it turns out to be simpler to first "lift" the problem to the level of a sample path large deviation problem, and then use the contraction principle to reduce to the original finite-dimensional problem. A "time" variable $x$ is introduced into the problem, where $\lfloor nx \rfloor$ balls have been thrown at time $x$, and $\Gamma_i^n(x)$ is equal to the fraction of cells containing $i$ balls at this time. We then follow a standard program: the large deviation properties of this Markov process are analyzed, the rate function $J$ on path space is obtained, and the rate function for the occupancy at time $\beta$ is then characterized as the solution to a variational problem involving $J$.

Although the program is standard, there are several very interesting features, both qualitative and technical, which distinguish this large deviation problem. We first describe some of the attractive qualitative features. Typically, one has a rate function on path space of the form $J(\phi) = \int L(\phi(x), \dot\phi(x))\, dx$, where the nonnegative function $L(\gamma, \xi)$ is jointly lower semi-continuous and convex in $\xi$ for each fixed $\gamma$. The large deviation properties of the process at time $\beta$ are then found by solving a variational problem of the form

$$\inf\{J(\phi) : \phi(\beta) = \omega\},$$

where $\omega$ is given and where there will also be constraints on the initial condition $\phi(0)$. In general, this problem does not have an explicit, closed form solution. One exception to this rule is the extraordinarily simple situation where $L(\gamma, \xi)$ does not depend on the state variable $\gamma$. In this case, Jensen's inequality implies that the minimizing trajectory is a straight line, and so the variational problem is actually finite dimensional. Another exception is the case of large deviation asymptotics of a small noise linear stochastic differential equation. In this case the variational problem takes the form of the classical linear quadratic regulator, and the explicit solution is well known from the theory of deterministic optimal control. However, in this case there is no need to "lift" the problem to the sample path level, since the distribution of the diffusion at any time is Gaussian with explicitly calculable mean and covariance. Other exceptions occur in large deviation problems from queuing theory, but in these problems the variational integrand is "sectionally" independent of $\gamma$, and one can show that the minimizing trajectories are the concatenation of a finite number of straight line segments (see, e.g., [1] and the references therein).

For the occupancy problem the variational integrand has a complicated state dependence [see (2.2)], reflecting the complicated dependence of the transition probabilities (in the process level version) on the state. Nonetheless, the rate function possesses a great deal of structure that can be heavily exploited. For example, the function $J$ turns out to be strictly convex on path space, and so all local minimizers in the variational problem are actually global minimizers. Perhaps even more surprising is the fact that explicit solutions to the Euler–Lagrange equations can be constructed, and as a consequence, the variational problem can be more-or-less solved explicitly (see the Appendix). Both of these properties follow from the fact that the variational integrand $L$ can be defined in terms of the famous *relative entropy* function (or divergence).

A technical novelty of the problem is the singular behavior of the transition probabilities of the underlying Markov process. Since only cells containing $j$ balls can become cells with $j+1$, it is clear that the transition probability corresponding to such an event scales linearly with $\Gamma_j^n(x)$, and in particular, that it vanishes at the boundary of the state space, when $\Gamma_j^n(x) = 0$. This poses no difficulty for the large deviations upper bound, but is an obstacle for the lower bound. (Some results that address lower bounds when rates go to zero include Chapter 8 of [26], which treats processes with "flat" boundaries, and recent general results in [27].) For the occupancy model, existing results provide a lower bound for open sets of trajectories that do not touch the boundary, because away from the boundaries the set $\{\xi : L(\Gamma, \xi) < \infty\}$ is independent of $\Gamma$. To deal with more general open sets we use a perturbation argument. We first show, using the strict convexity of $J$ and properties of the zero cost paths, that it is enough to prove large deviation lower bounds for open neighborhoods of trajectories that stay away from the boundary for all positive times. (Note that this still allows the trajectory to *start* on the boundary.) Loosely speaking, to prove the lower bound for sets of this type it is enough to show that given $a > 0$, there is $b > 0$ such that the probability that the process is at least distance $b$ from the boundary by time $b$ is bounded below by $\exp -na$. It turns out that these bounds can be easily established by exploiting an explicit representation for such probabilities that was proved in [11].

An outline of the paper is as follows. Section 2 states the main results of the paper. In the first part of Section 2 we construct the underlying stochastic process model, and state the corresponding sample path level LDP, as well as the LDP for the terminal distribution. The proof of the sample path LDP is given in the following section, although in Section 2 an additional heuristic argument for the form of the local rate function based on Sanov's theorem is provided. In the second part of Section 2, the terminal distribution rate function is characterized and explicit expressions for the sample path minimizers are presented. These constitute more-or-less complete solutions to the corresponding calculus of variations problem. The detailed proofs of these latter results are deferred to the Appendix.

Section 3 gives a proof of the sample path LDP. The solutions to the variational problem are exemplified in Section 4, where we show how specific questions regarding the occupancy problem can be answered. In particular, we work out the asymptotics for a number of examples, including an overflow problem and a partial coupon collection problem. Generalizations to our results are also described.

## 2. Main results.

2.1. *Statement of the LDP.* In this section we formulate the problem of interest and state the LDP. The proof will be given in Section 3. As noted in the Introduction, our focus is the asymptotic behavior of the occupancy problem. If $n$ denotes the total number of cells, then to have a nontrivial limit, the number of balls placed in the cells should scale linearly with $n$. We will place $\lfloor \beta n \rfloor$ balls in the cells, where $\beta \in (0, \infty)$ is a fixed parameter and $\lfloor a \rfloor$ denotes the integer part of $a$.

As will be seen in the sequel, the large deviation asymptotics of the occupancy should be treated by first lifting the problem to the level of sample path large deviations, and then using the contraction principle to reduce to the original problem. To do this, we introduce a "time" variable $x$ that ranges from 0 to $\beta$. At time $x$, one should imagine that $\lfloor nx \rfloor$ balls have been thrown. Thus, the occupancy process will be piecewise constant over intervals of the form $[i/n, i/n + 1/n)$. With this scaling of time, large deviation asymptotics can be obtained if we scale space by a factor of $1/n$. Thus, we define the *random occupancy process* $\Gamma^n(x) \doteq (\Gamma_0^n(x), \ldots, \Gamma_I^n(x), \Gamma_{I+}^n(x))$ by letting $\Gamma_i^n(x)$, $i = 0, \ldots, I$ denote $1/n$ times the number of cells with exactly $i$ balls at time $x$, and letting $\Gamma_{I+}^n(x)$ be $1/n$ times the number of cells with more than $I$ balls at time $x$. Note that $\Gamma^n$ takes values in the set of probability vectors on $I+2$ points: $S_I \doteq \{\gamma \in \mathbb{R}^{I+2} : \gamma_j \geq 0, 0 \leq j \leq I+1 \text{ and } \sum_{j=0}^{I+1} \gamma_j = 1\}$.

The process $\{\Gamma^n(i/n), i = 0, 1, \ldots\}$ is obviously Markovian. It will be convenient to work with the following "dynamical system" representation:

$$\Gamma^n\left(\frac{i+1}{n}\right) = \Gamma^n\left(\frac{i}{n}\right) + \frac{1}{n} b_i\left(\Gamma^n\left(\frac{i}{n}\right)\right),$$

where the independent and identically distributed random vector fields $\{b_i(\cdot), i = 0, 1, \ldots\}$ have distributions

$$P\{b_i(\gamma) = v\} = \begin{cases} \gamma_j, & \text{if } v = e_{j+1} - e_j, \ 0 \leq j \leq I, \\ \gamma_{I+1}, & \text{if } v = 0. \end{cases}$$

Here $e_j$ represents the vector in $\mathbb{R}^{I+2}$ whose $j$th element is unity and for which all other elements are zero. The occupancy after $\lfloor \beta n \rfloor$ balls have been thrown can then be represented by $n\Gamma^n(\beta)$.

As we have discussed, the large deviations behavior of $\Gamma^n(\beta)$ will be deduced from that of the process $\Gamma^n(\cdot)$. In order to state the large deviation asymptotics precisely, we should clarify the space in which $\Gamma^n(\cdot)$ takes values and the topology used on that space. As usual for processes with jumps of this sort, we use the Skorokhod space $D([0, \beta] : \mathbb{R}^{I+2})$ together with the Skorokhod topology [5], Chapter 4. However, the large deviations properties of $\{\Gamma^n\}$ are the same as for the processes $\{\bar{\Gamma}^n\}$, where $\bar{\Gamma}^n$ is defined to be the piecewise linear process which agrees with $\Gamma^n$ at times of the form $i/n$. For readers unfamiliar with the Skorokhod space and associated topology, the identical large deviations results hold for $\bar{\Gamma}^n$, save that the space of continuous functions and the sup norm topology are used instead.

To complete the statement of the LDP for $\{\Gamma^n\}$ we need some additional notation. A vector of (deterministic) *occupancy rates* $\theta(x) = (\theta_0(x), \ldots, \theta_I(x), \theta_{I+}(x))$ is a measurable mapping from $[0, \beta]$ to $S_I$. Intuitively these rates represent the rate at which balls flow into urns of a given occupancy level. Associated with each such vector of rates is the corresponding deterministic *occupancy function* $\gamma(x) = (\gamma_0(x), \ldots, \gamma_I(x), \gamma_{I+}(x))$, which is defined by the initial condition $\gamma(0)$ and the differential equations $\dot{\gamma}_0(x) = -\theta_0(x)$, $\dot{\gamma}_j(x) = \theta_{j-1}(x) - \theta_j(x)$ for $j = 1, \ldots, I$, and $\dot{\gamma}_{I+}(x) = \theta_I(x)$. These equations reflect the idea that the fraction of urns containing $i$ balls increases as balls enter $(i-1)$-occupied urns, and decreases as balls enter $i$-occupied urns. Defining the matrix

$$(2.1) \qquad M = \begin{pmatrix} -1 & 0 & 0 & \cdots & 0 \\ 1 & -1 & 0 & \cdots & 0 \\ 0 & 1 & -1 & \cdots & 0 \\ \vdots & \vdots & & \ddots & \\ 0 & \cdots & 1 & -1 & 0 \\ 0 & \cdots & 0 & 1 & 0 \end{pmatrix},$$

we can write $\dot{\gamma}(x) = M\theta(x)$ for all $x \in [0, \beta]$. Conversely, given a differentiable occupancy function $\gamma$, the corresponding rates are uniquely determined by $\dot{\gamma}(x) = M\theta(x)$ and the normalization $\sum_{i=0}^{I} \theta_i(x) + \theta_{I+}(x) = 1$. We will also be interested in the *cumulative occupancy function* $\psi$ with components $\psi_i(x) = \sum_{j=0}^{i} \gamma_j(x)$. Inspecting the cumulative sums of the rows of $M$ shows that $\dot{\psi}_i = -\theta_i$ for $i = 0, \ldots, I$. As more balls are thrown, the fraction of urns containing $i$ or fewer balls can only decrease, and the rate of decrease is the rate at which balls enter $i$-occupied urns. We will say that $\gamma$ is a *valid* occupancy function if $\gamma$ is absolutely continuous, if $\gamma(x)$ is a probability vector for all $x \in [0, \beta]$, and if its associated $\theta(x)$ is a probability vector for almost all $x \in [0, \beta]$. Note that the functions $\psi$, $\gamma$ and $\theta$ are interchangeable, in the sense that each can be derived from any of the others [given $\gamma(0)$ in the case of $\theta$]. Thus, we say that $\theta$ and $\psi$ are valid if the associated $\gamma$ is valid. The following lemma gives a direct characterization of validity in terms of $\psi$.

LEMMA 2.1. *A vector of $I+1$ continuous functions $\psi$, each of which maps $[0, \beta]$ to $[0, 1]$, is a valid cumulative occupancy path if and only if:*

(a) $\psi_i(x) \geq \psi_{i-1}(x)$,
(b) $\psi_i(x) \geq \psi_i(y)$,
(c) $\sum_{k=0}^{I}(\psi_k(x) - \psi_k(y)) \leq y - x$,

*for each $0 \leq i \leq I$ and $0 \leq x < y \leq \beta$.*

SKETCH OF THE PROOF. In the forward direction, the first condition plus the bound $\psi_I(x) \leq 1$ ensure that $\gamma(x)$ is a probability distribution. The second and third conditions together imply that $\psi$ is Lipschitz continuous with constant 1 and, hence, absolutely continuous. This implies the absolute continuity of $\gamma$. Since $\dot{\psi} = -\theta$, the third condition ensures that $\theta$ is almost always a probability distribution. The reverse direction proceeds similarly. □

Given $\gamma, \theta \in S_I$, let $D(\theta \| \gamma)$ denote relative entropy of $\theta$ with respect to $\gamma$. Thus,

$$D(\theta \| \gamma) \doteq \sum_{i=0}^{I+1} \theta_i \log(\theta_i / \gamma_i),$$

with the understanding that $\theta_i \log(\theta_i / \gamma_i) \doteq 0$ whenever $\theta_i = 0$, and that $\theta_i \log(\theta_i / \gamma_i) \doteq \infty$ if $\theta_i > 0$ and $\gamma_i = 0$. If $\gamma(x)$ is a valid occupancy function with corresponding occupancy rates $\theta(x)$, then we set

(2.2) $$J(\gamma) = \int_0^\beta D(\theta(x) \| \gamma(x)) \, dx.$$

In all other cases set $J(\gamma) = \infty$.

THEOREM 2.2. *Suppose that the sequence of initial conditions $\{\Gamma^n(0), n = 1, 2, \ldots\}$ is deterministic and that it converges to $\alpha \in S_I$ as $n \to \infty$. Then the sequence of processes $\{\Gamma^n, n = 1, 2, \ldots\}$ satisfies the LDP with rate function $J$. In other words, for any measurable set $A$ of trajectories, we have the large deviation lower bound*

$$\liminf_{n \to \infty} \frac{1}{n} \log P\{\Gamma^n \in A\} \geq -\inf\{J(\gamma) : \gamma \in A^\circ, \gamma(0) = \alpha\}$$

*and the large deviation upper bound*

$$\limsup \frac{1}{n} \log P\{\Gamma^n \in A\} \leq -\inf\{J(\gamma) : \gamma \in \bar{A}, \gamma(0) = \alpha\},$$

*and, moreover, for any compact set of initial conditions $K$ and $C < \infty$, the set*

$$\{\gamma : J(\gamma) \leq C, \gamma(0) \in K\}$$

*is compact.*

The proof of this result is provided in Section 3. However, a formal justification for the form of the rate function is as follows. Let $\delta > 0$ be a small time increment. Owing to the fact that $\Gamma^n$ varies slowly when compared to $\{b_i(\gamma)\}$, one expects

$$P\left(\frac{\Gamma^n(\delta) - \Gamma^n(0)}{\delta} \approx \eta \Big| \Gamma^n(0) = \gamma\right) \approx P\left(\frac{1}{n\delta} \sum_{j=0}^{\lfloor n\delta \rfloor - 1} b_j(\gamma) \approx \eta\right),$$

where the symbol "$\approx$" inside the probability means that the indicated quantities are within a small fixed constant $\varepsilon > 0$ of each other. Suppose that we interpret $\gamma$ as a probability measure on $\{e_0, \ldots, e_{I+1}\}$, and let $\{Y_j\}$ be independent and identically distributed (i.i.d.) with distribution $\gamma$. Then the sequence of i.i.d. random vectors $\{b_j(\gamma)\}$ can be realized by setting

$$b_j(\gamma) = \left\{\begin{matrix} e_{i+1} - e_i \\ 0 \end{matrix}\right\} \quad \Longleftrightarrow \quad Y_j = \left\{\begin{matrix} e_i \\ e_{I+1} \end{matrix}\right\},$$

that is,

$$b_j(\gamma) = MY_j.$$

By Sanov's theorem, for any probability vector $\theta \in S_I$,

$$P\left(\frac{1}{n\delta} \sum_{j=0}^{\lfloor n\delta \rfloor - 1} Y_j \approx \theta\right) \approx \exp -n\delta D(\theta \| \gamma).$$

Thus,

$$P\left(\frac{\Gamma^n(\delta) - \Gamma^n(0)}{\delta} \approx M\theta \Big| \Gamma^n(0) = \gamma\right) \approx \exp -n\delta D(\theta \| \gamma).$$

Approximating an arbitrary trajectory by a piecewise linear trajectory with nearly equal cost and using the Markov property, one expects

$$P(\Gamma^n \approx \gamma | \Gamma^n(0) = \alpha) \approx \exp-n \int_0^\beta D(\theta(t) \| \gamma(t)) \, dt,$$

where $\theta$ and $\gamma$ are related by $\dot\gamma = M\theta$, $\gamma(0) = \alpha$. Thus, the rate function on path space should be $J(\gamma)$.

The zero cost trajectories are of course the law of large numbers limits, and can easily be computed. For example, if $\alpha = (1, 0, \ldots)$ (so all cells are initially empty) and $i \leq I$, then $J(\gamma) = 0$ implies that $\gamma_i(t) = e^{-t} t^i / i!$. In other words, $\gamma(t)$ is the Poisson distribution with mean $t$, save that all mass corresponding to $i > I$ is collected together into the state $I+1$. Throughout this paper we will denote the Poisson distribution with mean $t$ by $\mathcal{P}(t)$, where $\mathcal{P}_i(t) = e^{-t} t^i / i!$.

We are primarily interested in the distribution of $\Gamma^n(\beta)$. The contraction principle (e.g., [11], Theorem 1.3.2) implies the following variational representation for the rate function for $\{\Gamma^n(\beta), n = 1, 2, \ldots\}$. Let $A(\alpha, \omega, \beta)$ denote the set of valid occupancy paths $\gamma$ on $[0, \beta]$ satisfying $\gamma(0) = \alpha$ and $\gamma(\beta) = \omega$.

COROLLARY 2.3. *Suppose that the sequence of initial conditions* $\{\Gamma^n(0), n = 1, 2, \ldots\}$ *is deterministic and that it converges to $\alpha$ as $n \to \infty$. Then the sequence of random vectors* $\{\Gamma^n(\beta), n = 1, 2, \ldots\}$ *satisfies the LDP with the convex rate function*

$$\mathcal{J}(\omega) = \inf\{J(\gamma) : \gamma \in A(\alpha, \omega, \beta)\}.$$

REMARK 2.1. Using the explicit formula for $\mathcal{J}(\omega)$ stated in Theorem 2.7 and the convexity of relative entropy in its first argument, it follows that $\mathcal{J}$ is, in fact, strictly convex.

REMARK 2.2. It is sometimes useful to show that the large deviation lower bound holds for sets with no interior relative to the ambient space. Such bounds can often be proved for processes, such as ours, that take values in a discrete lattice. An example would be a set $A \subset \{\gamma \in S_I : \gamma_{I+} = 0\}$, for which the minimizing trajectories must be polynomial extremals (defined after Theorem 2.6). Although we do not need such results in the present paper, it is worth observing that lower bounds of this kind can be proved.

2.2. *Characterization of the terminal rate function and minimizing paths.* The results of this section were obtained using calculus of variations techniques, the details of which are given in the Appendix. The presentation begins with the case in which all urns are initially empty because it appears in many applications and because it is a building block for the general case. In each case, first we give a characterization of $\mathcal{J}(\omega)$ as a minimal relative entropy, which may be computed explicitly in this form. We then give an explicit functional form for the sample path minimizer. The function form constain parameters determined by the solutions to fixed point equations.

Before stating these theorems, it is helpful to specify the domain over which $\mathcal{J}(\omega)$ is finite. We define an endpoint constraint to be a triple $(\alpha, \omega, \beta)$, where $\alpha, \omega \in S_I$ are the initial and terminal occupancies, and where $\beta > 0$ is the number of balls thrown per urn. We assume without loss that $\alpha_0 > 0$, that is, that some fraction of urns are initially empty, and we denote the index set of positive initial occupancies by $\mathcal{K} = \{k : \alpha_k > 0\}$. Since we do not distinguish between urns having more than $I$ balls, we can always suppose that no urns initially have more than $I+1$ balls, and denote the last element of $\alpha$ by $\alpha_{I+1}$. The last element of $\omega$ is denoted $\omega_{I+}$, signifying that it collects all urns with occupancy greater than or equal to $I+1$.

DEFINITION 2.1. A endpoint constraint $(\alpha, \omega, \beta)$ is feasible if the corresponding set of valid occupancy paths $A(\alpha, \omega, \beta)$ is nonempty.

LEMMA 2.4. An endpoint constraint $(\alpha, \omega, \beta)$ is feasible if and only if

$$\sum_{j=0}^{i} \alpha_j \geq \sum_{j=0}^{i} \omega_j, \qquad i = 0, \ldots, I \text{ (monotonicity)} \tag{2.3}$$

and

$$\sum_{i=0}^{I} i\omega_i + (I+1)\omega_{I+} \leq \sum_{k=0}^{I+1} k\alpha_k + \beta \qquad \text{(conservation)}. \tag{2.4}$$

This lemma is proved in the Appendix. Condition (2.3) relates to the fact that the $\psi_i(x)$ must decrease monotonically and (2.4) to a conservation constraint for the number of balls thrown. The right-hand side of the inequality equals the initial number of balls per urn, plus the additional balls per urn thrown up to time $\beta$, while the left-hand side is a lower bound on the number of balls per urn at time $\beta$.

In the treatment of general initial conditions, the following further definition will be useful.

DEFINITION 2.2. A set of feasible constraints $(\alpha, \omega, \beta)$ is *irreducible* if the monotonicity conditions (2.3) hold with strict inequality for all $i < I$. Otherwise, the constraints are termed reducible.

For a reducible set of constraints, let $i$ be the first index such that $\sum_{k=0}^{i} \alpha_k = \sum_{k=0}^{i} \omega_k$. This condition can only be met if no balls are ever thrown into $i$-occupied urns, and it follows that urns which initially contain $i$ balls or less will always contain $i$ balls or less. Thus, these urns may be treated in isolation from the urns containing more than $i$ balls. Furthermore, by subtracting $i+1$ balls from each urn in this latter set, we obtain another occupancy problem in standard form, that is, with $\alpha_0 > 0$. Continuing in this way, a given problem with constraints $(\alpha, \omega, \beta)$ may be divided into a finite number of isolated, irreducible subproblems. It is only necessary therefore to treat problems with irreducible constraints, and wherever necessary we suppose this to be the case.

2.2.1. *Empty initial conditions.* An important special case of the initial constraint is when all urns are initially empty, that is, $\alpha_0 = 1$ and $\alpha_i = 0$ for $i > 0$. We abuse notation and denote this case by $\alpha = 1$.

Define the set $F(1, \omega, \beta)$ to be the set of distributions $\pi$ on the nonnegative integers satisfying $\pi_i = \omega_i$ for $i = 0, \ldots, I$ and the constraint

$$\sum_{i=0}^{\infty} i \pi_i = \beta.$$

In the empty case, the conditions for feasibility of $(1, \omega, \beta)$ reduce to $\sum_{i=0}^{I} i \omega_i + (I+1)\omega_{I+} \leq \beta$, from which it follows that the set $F$ is nonempty if and only if $(1, \omega, \beta)$ is feasible. As we now discuss, the infimum of the rate function $J$ over paths $A(1, \omega, \beta)$ can be represented as an infimum of a relative entropy over distributions in $F(1, \omega, \beta)$.

THEOREM 2.5 (Terminal rate function, empty case). *Given the initial condition $\alpha = 1$, the rate function $\mathcal{J}(\omega) : S_I \to \mathbb{R}^+$ defined in Corollary 2.3 may be determined as*

$$\mathcal{J}(\omega) = \min_{\pi \in F(1,\omega,\beta)} D(\pi \| \mathcal{P}(\beta))$$

*if $(1, \omega, \beta)$ is feasible, and is otherwise infinite. The minimizing argument $\pi^* \in F(1, \omega, \beta)$ is unique.*

The above expression is of interest in its own right. Moreover, the optimal solution $\pi^*$ can be computed explicitly. Using Lagrange multipliers, one can show, in fact, that the solution takes the form $\pi_i^* = C\mathcal{P}_i(\rho\beta)$ for $i > I$ for some constants $C \geq 0$ and $\rho > 0$ that we refer to as twist parameters. Here $\rho$ is related to the Lagrange multiplier for the conservation constraint $\sum_{i=0}^{\infty} i\pi_i^* = \beta$, while $C$ is a normalization constant ensuring that $\sum_{i=0}^{\infty} \pi_i^* = 1$. These two constraints may be solved to determine $\rho$ and $C$. If we have the strict equality $\sum_{i=0}^{I} i\omega_i + (I+1)\omega_{I+} = \beta$, there are just enough balls to meet the terminal constraints and we may replace $I$ by $I+1$ if necessary to ensure that $\omega_{I+} = 1 - \sum_{i=0}^{I} \omega_i = 0$. In this case, it turns out that $F(1, \omega, \beta)$ has only one element; we then have $C = 0$, and we may take $\rho = 1$. Otherwise, we define $\rho$ to be the unique positive root of the equation

$$\text{(2.5)} \qquad \frac{\rho\beta - \sum_{i=0}^{I} i\mathcal{P}_i(\rho\beta)}{1 - \sum_{i=0}^{I} \mathcal{P}_i(\rho\beta)} = \frac{\beta - \sum_{i=0}^{I} i\omega_i}{1 - \sum_{i=0}^{I} \omega_i}.$$

Because of the strict inequality in the conservation condition (2.4), the right-hand side of the last equation is strictly greater than $I+1$. The left-hand side of the equation is the conditional mean $E[Y \mid Y > I]$ of a Poisson random variable $Y$ with mean $\rho\beta$. As a function of $\rho$, this conditional mean is a strictly monotonic and continuous map from $(0, \infty)$ to $(I+1, \infty)$ and, hence, the equation has exactly one positive root. $C$ is given by

$$\text{(2.6)} \qquad C = \frac{1 - \sum_{i=0}^{I} \omega_i}{1 - \sum_{i=0}^{I} \mathcal{P}_i(\rho\beta)} = \frac{\beta - \sum_{i=0}^{I} i\omega_i}{\rho\beta - \sum_{i=0}^{I} i\mathcal{P}_i(\rho\beta)}.$$

Evaluating the relative entropy of $\pi^*$ and $\mathcal{P}(\beta)$, the rate function of $\Gamma^n(\beta)$ can be expressed as

$$\text{(2.7)} \qquad \mathcal{J}(\omega) = \sum_{i=0}^{I} \omega_i \log \frac{\omega_i}{\mathcal{P}_i(\beta)} + \left(1 - \sum_{i=0}^{I} \omega_i\right)(\log C + (1-\rho)\beta) + \left(\beta - \sum_{i=0}^{I} i\omega_i\right)\log \rho.$$

The next theorem shows that the least cost path $\gamma^*$ satisfying $J(\gamma^*) = \mathcal{J}(\omega)$ can also be expressed explicitly in terms of the twist parameters $C$ and $\rho$. The least cost path is of interest for a number of reasons. First of all, the proof of Theorem 2.5 is obtained by evaluating (2.2), using the explicit form of $\gamma^*$. In a similar way, the least cost paths may be used as a tool in other problems of interest, such as for determining the rate function of $\Gamma^n(\beta)$ when $\beta$ is random. Second, the least cost path provides insight into the expected behavior of occupancy experiments, conditioned on the occurrence of a rare event. Third, they allow empirical estimation of rare event probabilities by change-of-measure importance sampling. A key step in proving the minimality of the proposed least cost path is to show that the path satisfies the Euler–Lagrange equations (defined in the Appendix).

THEOREM 2.6 (Globally minimizing path, empty case). *Suppose that $(1, \omega, \beta)$ are feasible constraints with empty initial conditions. The infimum of* $J(\gamma)$ over $A(1,\omega,\beta)$ is achieved on the occupancy path $\gamma \in A(1,\omega,\beta)$ defined by

$$\gamma_0(x) = Ce^{-\rho x} + \sum_{k=0}^{I}(\omega_k - C\mathcal{P}_k(\rho\beta))\left(1 - \frac{x}{\beta}\right)^k, \tag{2.8}$$

$$\gamma_i(x) = \frac{x^i}{i!}(-1)^i \gamma_0^{(i)}(x), \qquad 1 \leq i \leq I, \tag{2.9}$$

$$\gamma_{I+}(x) = 1 - \sum_{i=0}^{I} \gamma_i(x),$$

where $\rho > 0$ and $C \geq 0$ are twist parameters associated with the constraints. In addition, $\gamma$ satisfies the Euler–Lagrange equations.

Note that the entire path $\gamma(x)$ is completely determined by the empty component $\gamma_0(x)$ and its derivatives. In particular, the components $\gamma_i(x)$ are the terms in the Taylor expansion of $\gamma_0(x)$ about $x$, $\gamma_0(x+y) = \gamma_0(x) + y\gamma_0^{(1)}(x) + \cdots$ evaluated at time 0, that is, with $y = -x$. Note that $\gamma_0(x)$ is the sum of a polynomial and a single exponential term, so that the Taylor expansion always exists. When $C = 0$, there is no exponential term, and we say that $\gamma$ is a *polynomial* extremal. Otherwise, $C > 0$, we have an *exponential* extremal.

2.2.2. *General initial conditions.* Occupancy problems with general initial conditions may be thought of as a coupled set of problems with empty initial conditions. In particular, we consider the set of urns initially containing $k$ balls to form a class. The evolution of excess balls (beyond $k$) entering urns of this class may be denoted by occupancy functions of the form $\gamma_{k,j}(t)$ representing the fraction of balls initially having $k$ urns which have $k+j$ urns at time $t$. The fraction of urns containing $i$ balls in the overall system is obtained by summing contributions from all subproblem components with $k + j = i$.

As in the empty case, the rate function $\mathcal{J}(\omega)$ can be expressed as the solution to a minimization problem. Let $S_\infty$ denote the set of distributions on the nonnegative integers, and let $S_\infty^n$ denote the set of $n$-tuples of such distributions. Recall that $\mathcal{K}$ is set of indices $k$ such that $\alpha_k > 0$. We will denote an element of $S_\infty^{|\mathcal{K}|}$ by $\pi = \{\pi_{(0)}, \ldots, \pi_{(k)}, \ldots, \pi_{(K)}\}$, where for any $k \in \mathcal{K}$ a component distribution is denoted $\pi_{(k)} = \{\pi_{k,0}, \pi_{k,1}, \ldots\}$, and it is understood that the corresponding $\pi_{(k)}$ is omitted if $\alpha_k = 0$. Finally, let $F(\alpha, \beta, \omega)$ be the set of $\pi \in S_\infty^{|\mathcal{K}|}$ which satisfies the terminal constraints

$$\omega_i = \sum_{k \leq i, k \in \mathcal{K}} \alpha_k \pi_{k,i-k} \qquad \text{for all } 0 \leq i \leq I, \tag{2.10}$$

along with the conservation constraint

$$\sum_{k \in \mathcal{K}} \alpha_k \sum_{j=0}^{\infty} j \pi_{k,j} = \beta. \tag{2.11}$$

As in the empty case, it may be established that $F$ is nonempty if and only if $(\alpha, \omega, \beta)$ is feasible.

THEOREM 2.7 (Terminal rate function, general case). *The rate function $\mathcal{J}(\omega): S_I \to \mathbb{R}^+$ defined in Corollary 2.3 may be expressed*

$$\mathcal{J}(\omega) = \min_{\pi \in F(\alpha,\omega,\beta)} \sum_{k \in \mathcal{K}} \alpha_k D(\pi_{(k)} \| \mathcal{P}(\beta))$$

*whenever $(\alpha, \beta, \omega)$ are feasible, and is infinite otherwise. The minimizing argument $\pi^* \in F(\alpha, \omega, \beta)$ is unique.*

As discussed above, we may suppose the constraints are irreducible, and then, as shown in the Appendix, Lagrange multipliers will always exist for this problem. When we have strict inequality in the conservation condition (2.4) (the exponential case), the solution takes the form

$$\pi_{k,j}^* = \begin{cases} C_k \mathcal{P}_j(\rho\beta) W_{k+j}, & k \in \mathcal{K},\ k+j \leq I, \\ C_k \mathcal{P}_j(\rho\beta), & k \in \mathcal{K},\ k+j > I. \end{cases}$$

In the case of equality in (2.4) (the polynomial case), the corresponding form is

$$\pi_{k,j}^* = \begin{cases} D_k \mathcal{P}_j(\beta) W_{k+j}, & k \in \mathcal{K},\ j+k \leq I, \\ 0, & k+j > I. \end{cases}$$

As for empty initial conditions, $\rho$ may be associated with the conservation condition. The $W_i$ correspond to the terminal constraints $\omega_i$, and the $C_k, D_k$ are normalization constants. The constants $C_k$, $\rho$, $W_i$ and $D_i$ can all be computed numerically using Lagrangian methods for constrained optimization (see, e.g., [4]). Given these constants, the optimizing trajectory $\gamma(t)$ may be constructed explicitly. This construction may be most simply expressed in terms of the minimizing trajectories in the empty case. For $k \in \mathcal{K}$, denote the mean of $\pi_{(k)}^*$ by $\beta_k = \sum_{j=0}^{\infty} j \pi_{k,j}^*$, and define the terminal condition $\omega_{(k)} \in S_{I-k}$ by $\omega_{k,j} = \pi_{k,j}^*$, for $j = 0, \ldots, I-k$. Then the constraints $(1, \omega_{(k)}, \beta_k)$ are feasible constraints, and Theorem 2.6 determines the associated least cost paths, which we denote $\gamma_{(k)} = \{\gamma_{k,j}(x)\}$. The least cost paths for the subproblems combine to form the overall least cost path.

THEOREM 2.8 (Globally minimizing path, general case). *For irreducible feasible constraints $(\alpha, \omega, \beta)$, let $\pi^* \in S_\infty^{|\mathcal{K}|}$ be the unique minimizing distribution in Theorem 2.7, and let the functions $\gamma_{k,j}: [0, \beta_k] \to [0, 1]$ be the minimizing paths corresponding to the subproblems $(1, \omega_{(k)}, \beta_k)$. The infimum* *of $J(\gamma)$ over $A(\alpha,\omega,\beta)$ is achieved on the occupancy path $\gamma \in A(\alpha,\omega,\beta)$ defined by*

$$\gamma_i(x) = \sum_{k=0}^{i} \alpha_k \gamma_{k,i-k}(x\beta_k/\beta), \qquad i = 0, \ldots, I,$$

$$\gamma_{I+}(x) = 1 - \sum_{i=0}^{I} \gamma_i(x).$$

*In addition, $\gamma$ satisfies the Euler–Lagrange equations.*

**3. Proof of Theorem 2.2.** The purpose of this section is to prove the main large deviations result. We recall the processes and notation defined at the beginning of Section 2.

For $\zeta \in \mathbb{R}^{I+2}$ and $\gamma \in S_I$, define

$$\bar{H}(\gamma,\zeta) = \log(E[\exp(\zeta \cdot b_i(\gamma))]),$$

where "·" denotes inner product. Since the support of $b_i(\gamma)$ is bounded uniformly in $i$ and $\gamma$, there exists a function $h : \mathbb{R} \to \mathbb{R}$ such that $\bar{H}(\gamma,\zeta) \leq h(|\zeta|)$ for all $\zeta$ and $\gamma$. Also, since the distribution of $b_i(\gamma)$ is weakly continuous in $\gamma$, $\bar{H}(\gamma,\zeta)$ is jointly continuous. It follows from [10], Theorem 4.1, that the sequence $\{\Gamma^n, n = 1, 2, \ldots\}$ satisfies a large deviation upper bound with a rate function $\bar{J}$, which we now define. Let $\bar{L}$ be the Legendre–Fenchel transform of $\bar{H}(\gamma,\zeta)$ in $\zeta$:

$$\bar{L}(\gamma,\eta) = \sup_{\zeta \in \mathbb{R}^{I+2}} [\zeta \cdot \eta - \bar{H}(\gamma,\zeta)].$$

If $\gamma(x), 0 \leq x \leq \beta$ is an absolutely continuous function that takes values in $S_I$, then

$$\bar{J}(\gamma) = \int_0^\beta \bar{L}(\gamma(x),\dot{\gamma}(x)) \, dx.$$

If $\gamma$ is not absolutely continuous, then $\bar{J}(\gamma) = \infty$. ([10] assumes that the vector fields $b_i(\gamma)$ are defined for all $\gamma \in \mathbb{R}^{I+2}$. It is easy to check that we can extend the definition to this set with the bound $\bar{H}(\gamma,\zeta) \leq h(|\zeta|)$ and the continuity of $\bar{H}(\gamma,\zeta)$ preserved. However, if the process $\Gamma^n$ starts in $S_I$, then it stays in $S_I$, and so the exact form of the extension has no effect on the rate function.)

We recall the definition of the matrix $M$ given in (2.1). To complete the proof of the upper bound we must show that $\bar{J} = J$, where $J$ is defined as in (2.2). This will hold if we can show that $\bar{L}(\gamma,\eta)$ is finite only when $\eta = M\theta$ for a unique probability vector $\theta$, and that in this case

$$\bar{L}(\gamma,M\theta) = D(\theta\|\gamma).$$

It is well known that $\bar{L}(\gamma, \eta)$ is finite if and only if $\eta$ is in the convex hull of the support of $b_i(\gamma)$ (see, e.g., [11], Lemma 6.2.3(d)). Therefore if $\bar{L}(\gamma, \eta) < \infty$, then $\eta$ can be written as a convex combination of the form

$$\theta_0(e_1 - e_0) + \cdots + \theta_I(e_{I+1} - e_I),$$

where $\theta_j \geq 0$ and $\sum_{j=0}^{I} \theta_j \leq 1$. Since the vectors $\{(e_{j+1} - e_j), j = 0, 1, \ldots, I\}$ are linearly independent, these values are unique. Setting $\theta_{I+1} = 1 - \sum_{j=0}^{I} \theta_j$, we have $\eta = M\theta$ for a unique probability vector $\theta$. Now assume that $\eta$ takes this form. Then

$$\bar{L}(\gamma, M\theta) = \sup_{\zeta \in \mathbb{R}^{I+2}} [\zeta \cdot M\theta - \log(E[\exp(\zeta \cdot b_i(\gamma))])]$$

$$= \sup_{\zeta \in \mathbb{R}^{I+2}} \left[ M^T \zeta \cdot \theta - \log\left( \left[ \sum_{j=0}^{I} \gamma_j \exp(\zeta_{j+1} - \zeta_j) \right] + \gamma_{I+1} \right) \right]$$

$$= \sup_{\zeta \in \mathbb{R}^{I+2}} \left[ \sum_{j=0}^{I} (\zeta_{j+1} - \zeta_j)\theta_j - \log\left( \left[ \sum_{j=0}^{I} \gamma_j \exp(\zeta_{j+1} - \zeta_j) \right] + \gamma_{I+1} \right) \right].$$

Given any values $\mu_0, \ldots, \mu_{I+1}$, we can define $\zeta_0, \ldots, \zeta_{I+1}$ recursively by $\zeta_0 = -\mu_{I+1}$ and $\zeta_{j+1} - \zeta_j = \mu_j - \mu_{I+1}$. With these definitions, it is apparent that the last display is equal to

$$\sup_{\mu \in \mathbb{R}^{I+2}} \left[ \sum_{j=0}^{I} \mu_j \theta_j - \mu_{I+1}(1 - \theta_{I+1}) - \log\left( \left[ \sum_{j=0}^{I} \gamma_j \exp(\mu_j - \mu_{I+1}) \right] + \gamma_{I+1} \right) \right]$$

$$= \sup_{\mu \in \mathbb{R}^{I+2}} \left[ \sum_{j=0}^{I} \mu_j \theta_j - \mu_{I+1}(1 - \theta_{I+1}) + \mu_{I+1} - \log\left( \sum_{j=0}^{I+1} \gamma_j \exp \mu_j \right) \right]$$

$$= \sup_{\mu \in \mathbb{R}^{I+2}} \left[ \sum_{j=0}^{I+1} \mu_j \theta_j - \log\left( \sum_{j=0}^{I+1} \gamma_j \exp \mu_j \right) \right].$$

According to the Donsker–Varadhan variational formula for relative entropy (e.g., [11], Lemma 1.4.3(a)), the last display equals $D(\theta \| \gamma)$, thereby completing the proof of the upper bound.

We turn now to the proof of the lower bound. In this proof we will assume that $\gamma_0(0) > 0$. Since $\gamma_0$ is always nonincreasing, $\gamma_0(0) = 0$ implies that $\gamma_0(x) = 0$ for all $x$ and, therefore, under this condition the first component plays no significant role. A proof analogous to the one given below applies when $\gamma_0(0) = 0$, where the role of $\gamma_0(0)$ here is played by the first positive component of $\gamma(0)$.

Let $P_{\Gamma^n(0)}$ [resp. $E_{\Gamma^n(0)}$] denote probability (resp. expected value) given a deterministic initial occupancy $\Gamma^n(0)$. To prove the large deviation lower bound, it suffices to show that given any $\varepsilon > 0$ and $\delta > 0$, there is $\eta > 0$ such that for any initial occupancies satisfying $|\Gamma^n(0) - \gamma(0)| < \eta$,

$$(3.1) \quad \liminf_{n \to \infty} \frac{1}{n} \log P_{\Gamma^n(0)} \left( \sup_{0 \leq x \leq \beta} |\Gamma^n(x) - \gamma(x)| < \delta \right) \geq -J(\gamma) - \varepsilon.$$

Of course this inequality is trivial if $J(\gamma) = \infty$, and so we assume that $J(\gamma) < \infty$.

As we have remarked, a source of difficulty is the singular behavior of the transition rates of the process when $\Gamma^n$ is near the boundary of $S_I$. We first show that this can be avoided at all times save $x = 0$. To do this, we show that for any $a > 0$, there exist $b > 0$, $K \in \mathbb{N}$ and an occupancy function $y$ such that $y(0) = \gamma(0)$, $\sup_{0 \leq x \leq \beta} |y(x) - \gamma(x)| < a$, $y_j(x) > bx^K$ for all $j = 0, 1, \ldots, I, I+$ and $0 < x \leq \beta$, and such that

$$J(y) \leq J(\gamma).$$

Consider the zero cost trajectory defined by

$$\dot{z}(x) = Mz(x), \qquad z(0) = \gamma(0).$$

We have the following expression for $z_j$ when $j \leq I$:

$$(3.2) \quad z_j(x) = \left[ \sum_{k=0}^{j} \gamma_k(0) x^{(j-k)}/(j-k)! \right] e^{-x}.$$

It is easy to check from this explicit formula that $z_j(x) > \bar{b} x^K$ for some $\bar{b} > 0$, $K = I$, and all $j = 0, 1, \ldots, I, I+$ and $0 < x \leq \beta$. For $\rho \in (0, 1)$, let $y^\rho = \rho z + (1 - \rho)\gamma$. Then $y^\rho$ is the occupancy function that corresponds to the rate $\rho z + (1 - \rho)\theta$. Using the joint convexity of relative entropy in both variables ([11], Lemma 1.4.3(b)) and the fact that $D(z\|z) = 0$, we have

$$J(y^\rho) = \int_0^\beta D(\rho z(x) + (1-\rho)\theta(x) \| \rho z(x) + (1-\rho)\gamma(x)) \, dx$$
$$\leq \rho \int_0^\beta D(z(x)\|z(x)) \, dx + (1-\rho) \int_0^\beta D(\theta(x)\|\gamma(x)) \, dx$$
$$= (1-\rho) J(\gamma).$$

All required properties are then obtained by letting $y = y^\rho$ for suitably small $\rho \in (0, 1)$.

It follows that in proving the lower bound, we can assume without loss of generality that for some fixed constants $b > 0$ and $K \in \mathbb{N}$, $\gamma_i(x) \geq bx^K$ for all $x \in [0, \beta]$. We now return to the proof of the lower bound. Our first objective is to show that the process can be moved into a small neighborhood of $\gamma(\tau)$ (for $\tau > 0$ small) with sufficiently high probability. Given $\tau \in (0, \beta], \sigma > 0$, and $\varepsilon > 0$, define

$$h(\gamma) = \begin{cases} 0, & |\gamma - \gamma(\tau)| < \sigma/2, \\ 2\varepsilon, & \text{else.} \end{cases}$$

For $n$ large enough that $|\gamma(\lfloor n\tau \rfloor/n) - \gamma(\tau)| \leq \sigma/2$ (and independent of $\tau \in (0, \beta])$, we have the inequality

$$P_{\Gamma^n(0)}(|\Gamma^n(\lfloor n\tau \rfloor/n) - \gamma(\lfloor n\tau \rfloor/n)| < \sigma) + e^{-2n\varepsilon}$$
$$\geq P_{\Gamma^n(0)}(|\Gamma^n(\lfloor n\tau \rfloor/n) - \gamma(\tau)| < \sigma/2) + e^{-2n\varepsilon}$$
$$\geq E_{\Gamma^n(0)}(\exp -nh(\Gamma^n(\lfloor n\tau \rfloor/n))).$$

We next exploit a representation for exponential integrals that will give us an explicit lower bound on the last quantity. Consider a process $\bar{\Gamma}^n(x)$ constructed as follows. The process dynamics are of the same general structure as those of $\Gamma^n$, save that $b_i(\Gamma^n(i/n))$ is replaced by a sequence $\bar{b}^n_i$:

$$\bar{\Gamma}^n\left(\frac{i+1}{n}\right) = \bar{\Gamma}^n\left(\frac{i}{n}\right) + \frac{1}{n}\bar{b}^n_i, \qquad \bar{\Gamma}^n(0) = \Gamma^n(0).$$

Furthermore, the distribution of $\bar{b}^n_i$ is allowed to depend in any measurable way upon the set of values $\{\bar{\Gamma}^n(j/n), 0 \leq j \leq i\}$. Let $\mu_\gamma$ denote the distribution of $b_i(\gamma)$, and (without explicitly exhibiting all the dependencies) let $\bar{\mu}^n_i$ denote the (random) distribution of $\bar{b}^n_i$, given $\{\bar{\Gamma}^n(j/n), 0 \leq j \leq i\}$. We let $\bar{E}_{\Gamma^n(0)}$ denote expectation on the space that supports these processes. It follows from [11], Theorem 4.3.1, that

$$-\frac{1}{n}\log E_{\Gamma^n(0)}\left(\exp -nh\left(\Gamma^n\left(\frac{\lfloor n\tau \rfloor}{n}\right)\right)\right)$$
$$= \inf \bar{E}_{\Gamma^n(0)}\left[h\left(\bar{\Gamma}^n\left(\frac{\lfloor n\tau \rfloor}{n}\right)\right) + \frac{1}{n}\sum_{i=1}^{\lfloor n\tau \rfloor -1} D(\bar{\mu}^n_i \| \mu_{\bar{\Gamma}^n(i/n)})\right],$$

where the infimum is over all such processes $\bar{\Gamma}^n(x)$. In order to obtain a lower bound, we now simply insert a particular choice for the random variables $\bar{b}^n_i$. We can write $\gamma(\tau) - \gamma(0) = Mv\tau$ for some probability vector $v$. Define a process $\bar{b}^n_i$ as follows:

$$\bar{b}^n_i = \begin{cases} e_j - e_{j-1}, & \text{if } \left\lfloor \sum_{k=0}^{j-1} v_k \tau n \right\rfloor \leq i \leq \left\lfloor \sum_{k=0}^{j} v_k \tau n \right\rfloor - 1, \text{ for } 0 \leq j \leq I, \\ 0, & \text{if } \left\lfloor \sum_{k=0}^{I} v_k \tau n \right\rfloor \leq i \leq \lfloor \tau n \rfloor - 1. \end{cases}$$

In other words, $\bar{b}^n_i$ defines a deterministic, discrete time approximation to the continuous time occupancy rate process that uses $e_0 - e_1$ for an amount of time $v_0\tau$, $e_1 - e_2$ for an amount of time $v_1\tau$ and so on. This continuous time process will move the occupancy process from $\gamma(0)$ to $\gamma(\tau)$ at time $\tau$. If $e_j - e_{j-1}$ is used at the discrete time step $i$, then since $\bar{\mu}_i^n$ concentrates its mass on $e_j - e_{j-1}$, the cost is

$$(3.3) \quad D(\bar{\mu}_i^n \| \mu_{\bar{\Gamma}^n(i/n)}) = \log\left(\frac{1}{\bar{\Gamma}_{j-1}^n(i/n)}\right) = -\log\left(\bar{\Gamma}_{j-1}^n\left(\frac{i}{n}\right)\right).$$

The process $\bar{\Gamma}^n(i/n)$ possesses important monotonicity and convergence properties. Since $\bar{\Gamma}_{j-1}^n((i+1)/n) - \bar{\Gamma}_{j-1}^n(i/n) = -1/n$ for $\lfloor \sum_{k=0}^{j-1} v_k \tau n \rfloor \leq i \leq \lfloor \sum_{k=0}^{j} v_k \tau n \rfloor - 1$,

$$\bar{\Gamma}_{j-1}^n\left(\frac{i}{n}\right) \downarrow \bar{\Gamma}_{j-1}^n\left(\frac{1}{n}\left\lfloor \sum_{k=0}^{j} v_k \tau n \right\rfloor\right)$$

as $i \uparrow \lfloor \sum_{k=0}^{j} v_k \tau n \rfloor$. In addition, because the $(j-1)$st component is never modified when $i \geq \lfloor \sum_{k=0}^{j} v_k \tau n \rfloor$, it follows that

$$\bar{\Gamma}_{j-1}^n\left(\frac{1}{n}\left\lfloor \sum_{k=0}^{j} v_k \tau n \right\rfloor\right) \to \gamma_{j-1}(\tau)$$

as $n \to \infty$ and $\eta \to 0$. Furthermore, as observed previously, (3.2) implies the existence of $b > 0$ and $K \in \mathbb{N}$ such that $\gamma_j(\tau) \geq b\tau^K$ for $j = 0, \ldots I$. Thus, at any given time step $i$ we have a strictly positive lower bound on the relevant component of $\bar{\Gamma}^n$, which in turn provides a strictly finite upper bound on the corresponding relative entropy cost. Indeed, it follows from (3.3) and $\gamma_j(\tau) \geq b\tau^K$ that for all sufficiently large $n$ and small $\eta > 0$, there are $C_1, C_2 < \infty$ (and independent of $\tau$) such that whenever $|\Gamma^n(0) - \gamma(0)| < \eta$, for all $i$,

$$D(\bar{\mu}_i^n \| \mu_{\bar{\Gamma}^n(i/n)}) \leq C_1[-\log \tau^K] \leq -C_2 \log \tau.$$

In addition, as $n \to \infty$ and $\eta \to 0$,

$$\bar{\Gamma}^n\left(\frac{1}{n}\lfloor \tau n \rfloor\right) \to \gamma(\tau).$$

By the Lebesgue dominated convergence theorem, for all sufficiently large $n$ and sufficiently small $\eta > 0$, $|\Gamma^n(0) - \gamma(0)| < \eta$ implies

$$-\frac{1}{n} \log E_{\Gamma^n(0)}\left(\exp -nh\left(\Gamma^n\left(\frac{\lfloor n\tau \rfloor}{n}\right)\right)\right) \leq -C_2 \tau \log \tau.$$

We now choose $\tau > 0$ so that $-C_2 \tau \log \tau \leq \varepsilon/2$. Choosing $\tau > 0$ smaller if need be, we can also guarantee that $|\Gamma^n(x) - \gamma(x)| \leq \delta$ for all $x \in [0, \tau]$ w.p.1 if $|\Gamma^n(0) - \gamma(0)| < \eta$ and $\eta > 0$ is sufficiently small. The following bound is

therefore valid for the given $\tau > 0$: for any $\sigma > 0$ and all sufficiently small $\eta > 0$, $|\Gamma^n(0) - \gamma(0)| < \eta$ implies

$$\liminf_{n \to \infty} \frac{1}{n} \log P_{\Gamma^n(0)}\left(\left|\Gamma^n\left(\frac{\lfloor n\tau \rfloor}{n}\right) - \gamma\left(\frac{\lfloor n\tau \rfloor}{n}\right)\right| < \sigma,\right.$$

$$\left. \sup_{x \in [0,\tau]} |\Gamma^n(x) - \gamma(x)| < \delta \right) \geq -\frac{\varepsilon}{2}.$$

Note that the asymptotic lower bound on the normalized log of the probability is independent of $\sigma > 0$. To obtain the lower bound for all $x \in [0, \beta]$, we will use the Markov property and an existing lower bound for paths which avoid that boundary. This latter lower bound will hold uniformly in a neighborhood of the initial condition $\gamma(\tau)$. Since we do not know a priori how small this neighborhood must be, it is important that the lower bound in the last display should be independent of $\sigma > 0$.

Now choose $\zeta \in (0, \delta]$ such that $\gamma(x)$ is at least distance $2\zeta$ from the boundary of $S_I$ for all $x \in [\tau, \beta]$. Recall that when considered as a function of $\gamma$, the distribution of $b_i(\gamma)$ is continuous in the weak topology, and moreover that the support of this distribution is independent of $\gamma$ so long as $\gamma_i > 0$ for all $i \in \{0, 1, \ldots, I+1\}$ (i.e., $\gamma \in S_I^\circ$, where $S_I^\circ$ denotes interior relative to the smallest affine space that contains $S_I$). It then follows from Proposition 6.6.1 of [11] (see also the discussion on [11], page 165, regarding uniformity) that $\{\Gamma^n, n = 1, 2, \ldots\}$ satisfies the following uniform large deviations lower bound: given any $\varepsilon > 0$ and $\zeta > 0$ defined above, there is $\sigma > 0$ such that as long as $|\Gamma^n(\lfloor n\tau \rfloor/n) - \gamma(\lfloor n\tau \rfloor/n)| < \sigma$,

$$\liminf_{n \to \infty} \frac{1}{n} \log P_{\Gamma^n(\lfloor n\tau \rfloor/n), \lfloor n\tau \rfloor/n}\left( \sup_{\tau \leq x \leq \beta} |\Gamma^n(x) - \gamma(x)| < \zeta \right) \geq -J(\gamma) - \frac{\varepsilon}{2},$$

where $P_{\Gamma^n(\lfloor n\tau \rfloor/n), \lfloor n\tau \rfloor/n}$ denotes probability given the occupancy levels $\Gamma^n(\lfloor n\tau \rfloor/n)$ at time $\lfloor n\tau \rfloor/n$. Proposition 6.6.1 of [11] assumes Condition 6.3.2. It is worth noting that in the present setting this condition holds with the particularly simple choice $\tilde{\beta} = \gamma$ (using the notation of [11]).

The lower bound (3.1) now follows by the Markov property and the last two displays. The proof that $J$ has compact level sets is as in [10], and therefore omitted.

**4. Examples and extensions.** In this section we apply the results of the previous sections to three different occupancy problems. We show how the parameters of interest may be computed numerically and plot the solutions to the associated calculus of variations problems. In Section 4.4 we list some other asymptotic problems of interest which can be solved by relatively straightforward generalizations of the results presented in this paper.

In the calculus of variations problems solved in Section 2, precise initial and terminal points were always given. In typical applications, one is interested in the minimum value of the rate function over a constraint set. When the constraint set is sufficiently simple, such problems may still be solved easily using the tools provided in Section 2.2. The three problems of this section are of this type.

Suppose that the event of interest is that the random endpoint $\Gamma^n(\beta)$ should lie in a terminal constraint set $\Omega$. To apply the LDP for $\{\Gamma^n(\beta)\}$ established by Corollary 2.3, one must compute exponents of the form

$$\mathcal{J}(\Omega) = \inf_{\omega \in \Omega} \mathcal{J}(\omega).$$

Using Theorem 2.7, we can write

(4.1)
$$\mathcal{J}(\Omega) = \inf_{\omega \in \Omega} \inf_{\pi \in F(\alpha,\omega,\beta)} \sum_{k=0}^{K} \alpha_k D(\pi_{(k)} \| \mathcal{P}(\beta))$$
$$= \inf_{\pi \in F(\alpha,\Omega,\beta)} \sum_{k=0}^{K} \alpha_k D(\pi_{(k)} \| \mathcal{P}(\beta)),$$

where we have abused notation to define $F(\alpha, \Omega, \beta) = \bigcup_{\omega \in \Omega} F(\alpha, \omega, \beta)$.

In many cases, for example, when the terminal set $\Omega$ is convex and defined by linear constraints, the exponent $\mathcal{J}$ can be computed directly from (4.1) using Lagrange multipliers. That is, one solves minimization problems of the type given in Theorem 2.7, but with the endpoint constraints (2.10) replaced by constraints defining $\Omega$. This is the approach used in the second and third example below. We do not prove that appropriate Lagrange multipliers always exist; if needed, existence may be established using methods similar to those used in the Appendix for the case $\Omega = \{\omega\}$. Because of the convexity of relative entropy in (4.1), a local minimum is always a global minimum over convex sets. Hence, in any particular scenario with convex $\Omega$, it is sufficient to establish a local minimum by numerically computing a set of Lagrange multipliers.

An alternative approach for computing $\mathcal{J}(\Omega)$ could be to return to the sample path level and use *natural boundary conditions* on the extremal curves (see [25]).

A set $\Omega$ with interior $\Omega^\circ$ and closure $\bar{\Omega}$ is a $\mathcal{J}$-continuity set if and only if $\inf_{\omega \in \Omega^\circ} \mathcal{J}(\omega) = \inf_{\omega \in \bar{\Omega}} \mathcal{J}(\omega)$. For such sets the large deviations lower and upper bounds coincide so that $-\lim 1/n \log P(\Gamma^n(\beta) \in \Omega) = \inf_{\omega \in \Omega} \mathcal{J}(\omega)$. It may readily be verified in each of the following examples that the event of interest is a $\mathcal{J}$-continuity set.

4.1. *The classical occupancy problem.* In the classical occupancy problem, the urns are initially empty and one only distinguishes between empty and occupied urns, or in other words, $I = 0$. The associated large deviations problem was solved using a sample path approach in [28]; we show here how this case may be obtained with our results. We might be interested in the probability of having an unusually large number of empty urns $\Gamma_0^n(\beta) > \omega_0 > e^{-\beta}$, or an unusually small number of empty urns $\Gamma_0^n(\beta) < \omega_0 < e^{-\beta}$. In either case, the calculus of variation problem is that in which $\gamma_0(\beta) = \omega_0$. From (2.6), it is immediate that $C = 1/\rho$ and

$$\gamma_0(x) = \frac{1}{\rho} e^{-\rho x} + 1 - \frac{1}{\rho}.$$

Using (2.5), we find that $\rho$ is determined by the unique nonnegative solution to

$$\rho(1 - \omega_0) = 1 - e^{-\beta\rho}.$$

Finally, (2.7) provides a simple expression for $\mathcal{J}(\omega)$ in terms of $\omega_0$, $\beta$, $C$ and $\rho$.

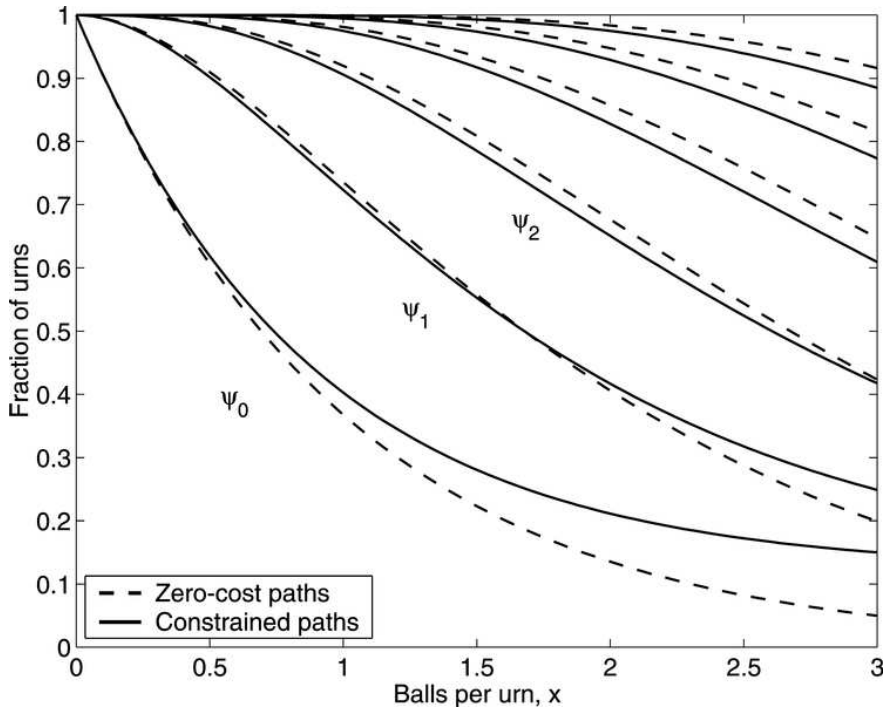

FIG. 1. *Cumulative urn occupancies $\psi_0$ through $\psi_5$ for the classical occupancy problem, including unconstrained paths (dashed curves) and constrained paths (solid curves).*

In addition, our analysis gives the sample paths for the higher occupancies, conditioned on an unusual number of empty urns, namely, $\gamma_i(x) = \mathcal{P}_i(\rho x)/\rho$ for $i > 0$. Figure 1 depicts the first five occupancy levels as a function of balls per urn thrown $x$. In this example, the terminal fraction of empty urns $\omega_0 = 0.15$ is three times larger than the expected value $\mathcal{P}_0(3.0) = 0.05$.

4.2. *The overflow problem.* In the overflow problem an urn is considered full when it reaches a finite capacity $I > 0$. Once an urn has been filled, successive balls thrown to that urn fall to the floor. For specificity, we consider the problem of determining the probability that an unusually large number of balls end up on the floor, and assume the empty initial conditions $\alpha_0 = 1$.

The problem can be handled using urns of infinite capacity in the following way. The number of balls that would have fallen on the floor in a finite capacity system is the number of balls in urns with occupancy greater than $I$, minus $I$ times the number of such urns. When $r$ balls have been thrown, the random number of overflowing balls $w(r)$ is thus

$$w(r) = r - n\Gamma_1(r/n) - 2n\Gamma_2(r/n) - \cdots - In\Gamma_I(r/n) - In\Gamma_{I+}(r/n)$$

or, since $\Gamma_{I+} = 1 - \sum_{i=0}^{I} \Gamma_i$,

$$w(r) = r - nI + n\sum_{i=0}^{I}(I-i)\Gamma_i(r/n).$$

In order to compute

$$J_O(\eta, \beta) \doteq -\lim_n \frac{1}{n} \log P\left(\frac{w(\beta n)}{n} > \eta\right),$$

we therefore consider sample paths which satisfy the end constraint

(4.2) $$\sum_{i=0}^{I}(I-i)\gamma_i(\beta) \geq \eta + I - \beta \doteq \zeta.$$

Note that $\zeta$ can be interpreted as the average spare capacity per urn, which must satisfy the bounds $[I - \beta]^+ \leq \zeta < I$, and that the average overflow satisfies $[\beta - I]^+ \leq \eta < \beta$. Assuming that $\eta$ (and $\zeta$) is larger than would be expected in the zero cost case, the minimum large deviations exponent will be achieved with equality in the constraint (4.2). Assuming that $\eta > 0$, we are in the exponential case, and the large deviations exponent $J_O(\eta, \beta)$ will be given by minimizing the divergence between $\pi$ and $\mathcal{P}(\beta)$, under the linear equality constraints

(4.3) $$\sum_{i=0}^{\infty} \pi_i = 1, \quad \sum_{i=0}^{\infty} i\pi_i = \beta \quad \text{and} \quad \sum_{i=0}^{I}(I-i)\pi_i = \zeta.$$

We introduce the Lagrangian

$$\mathcal{L}(\pi; y, z, \lambda) = \sum_{i=0}^{\infty} \pi_i \log \frac{\pi_i}{\mathcal{P}_i(\beta)} + y\left(1 - \sum_{i=0}^{\infty} \pi_i\right)$$
$$+ z\left(\beta - \sum_{i=0}^{\infty} i\pi_i\right) + \lambda\left(\zeta - \sum_{i=0}^{I}(I-i)\pi_i\right)$$

with Lagrange multipliers $y$, $z$ and $\lambda$. On differentiating, it follows that the minimizing distribution $\pi^*$ should satisfy the conditions

$$\log \pi_i^* = \log \mathcal{P}_i(\beta) - 1 + y + iz + (I-i)^+ \lambda.$$

In terms of variables $C, \rho$ and $\nu$, we may write

$$\pi_i^* = C\mathcal{P}_i(\rho\beta) \qquad \text{for } i \geq I$$

and

$$\pi_i^* = C\mathcal{P}_i(\rho\beta)\nu^{I-i} \qquad \text{for } i \leq I.$$

The distribution is conditionally Poisson $\rho\beta/\nu$ for $i \leq I$ and conditionally Poisson $\rho\beta$ for $i \geq I$. For convenience, we introduce the notation

$$Q_I(\rho) \doteq \sum_{i=I}^{\infty} \mathcal{P}_i(\rho\beta) = \frac{1}{\rho\beta} \sum_{i=I+1}^{\infty} i\mathcal{P}_i(\rho\beta),$$

$$R_I(\rho, \nu) \doteq \nu^I e^{-\rho\beta(1-1/\nu)} \sum_{i=0}^{I-1} \mathcal{P}_i\left(\frac{\rho\beta}{\nu}\right)$$

$$= \frac{\nu^{I+1}}{\rho\beta} e^{-\rho\beta(1-1/\nu)} \sum_{i=1}^{I} i\mathcal{P}_i\left(\frac{\rho\beta}{\nu}\right).$$

Since $\pi^*$ must satisfy the three linear constraints (4.3), the constants $C$, $\rho$ and $\nu$ must solve the equations

$$C(R_I(\rho, \nu) + Q_I(\rho)) = 1,$$
$$C\left(\frac{\rho\beta}{\nu}R_I(\rho, \nu) + \rho\beta Q_I(\rho)\right) = \beta,$$
$$C\left(I\mathcal{P}_I(\rho\beta) + \left(I - \frac{\rho\beta}{\nu}\right)R_I(\rho, \nu)\right) = \zeta.$$

There can be at most one positive triple $(C, \rho, \nu)$ satisfying these equations, since each such triple identifies a local minimum of $D(\pi \| \mathcal{P}(\beta))$ for $\pi$ in a convex set, and there can be only one such minimum. The equations can be solved numerically in a number of ways to obtain $C$, $\rho$ and $\nu$. For example, from the first constraint, $C$ can be expressed $C = (R_I + Q_I)^{-1}$. Substituting this expression into the second and third constraints, we obtain the equations

$$(\rho/\nu - 1)R_I(\rho,\nu) + (\rho - 1)Q_I(\rho) = 0,$$
$$(I - \rho\beta/\nu - \zeta)R_I(\rho,\nu) - \zeta Q_I(\rho) + I\mathcal{P}_I(\rho\beta) = 0.$$

Each equation implicitly defines a curve $\nu$ as a function of $\rho$, and the intersection of the two curves gives the desired $(\rho, \nu)$. We note that larger than expected values of $\zeta$ will lead to $\nu > \rho > 1$, while smaller than expected values of $\zeta$ give $\nu < \rho < 1$.

The large deviations exponent for the overflow problem may then be expressed

$$J_O(\eta, \beta) = D(\pi^* \| \mathcal{P}(\beta))$$
$$= \sum_{i=0}^{\infty} \pi_i^* \log(Ce^{\beta - \rho\beta} \rho^i \nu^{(I-i)^+})$$
$$= \log C + \beta(1 - \rho) + \beta \log \rho + \zeta \log \nu.$$

4.3. *Partial coupon collection, with initial conditions.* In the coupon collector's problem, the urns represent the $n$ types of coupons that are required to form a complete collection. The placement of a ball in a given urn corresponds to choosing a new coupon at random, and the problem is to see how many coupons must be collected before $I + 1$ complete sets are obtained. This event corresponds to the constraint $\omega_i = 0, i \leq I$.

In this section we solve a generalization of this problem. Beginning from nonempty initial conditions (a collection already in progress), we collect $\beta n$ additional coupons with the goal of obtaining more than $I$ coupons of as many types as possible. We want to determine how likely it is that number of types for which we have collected $I$ or fewer coupons is less than $\xi n$.

In terms of the urn problems we have considered, we are given initial occupancies $\alpha$ and wish to compute

$$J_C(\alpha, \beta, \xi) \doteq -\lim_n \frac{1}{n} \log P\left(\sum_{i=0}^{I} \Gamma_i^n(\beta) < \xi\right),$$

where

$$\xi < \sum_{k=0}^{K} \alpha_k \sum_{i=0}^{I-k} \mathcal{P}_i(\beta)$$

is an unusually small number of low occupancy urns.

The exponent $J_C$ will be given by computing $\mathcal{J}(\omega)$ as defined in Theorem 2.7 subject to the conservation constraint (2.11), replacing the terminal conditions (2.10) with the single constraint

$$\sum_{k=0}^{K} \alpha_k \sum_{j=0}^{I-k} \pi_{k,j} = \xi,$$

where $K \leq I$ since any sets which are initially complete may be left out of the problem. After constructing a Lagrangian and differentiating, we find that the minimizing solution must be of the form

$$\omega_{k,j} = \pi_{k,j}^* = \begin{cases} C_k \mathcal{P}_j(\rho\beta)W, & j+k \leq I, \\ C_k \mathcal{P}_j(\rho\beta), & j+k > I. \end{cases}$$

As in the previous example, the unknown constants may be determined by substituting the given form of the solution into the constraint equations, and solving the resulting system of equations. In terms of these constants, the large deviations exponent may be expressed

$$J_C(\alpha, \beta, \xi) = \beta(1 - \rho + \log \rho) + \xi \log W + \sum_{k=0}^{K} \alpha_k \log C_k.$$

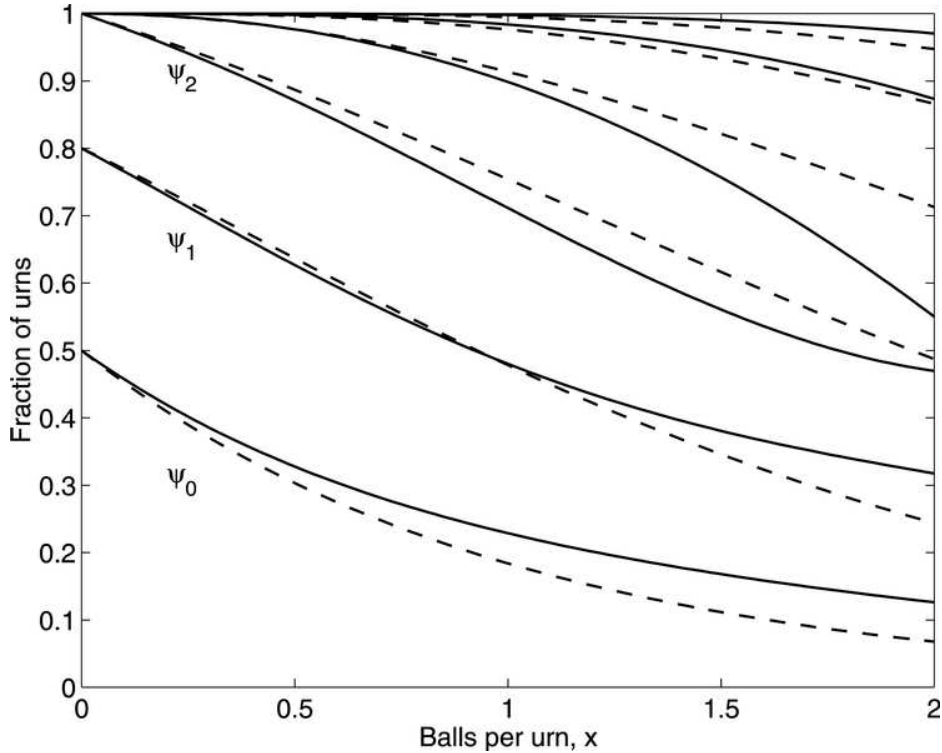

FIG. 2. *Cumulative urn occupancies $\psi_0$ through $\psi_5$, for a partial coupon collectors problem, including unconstrained paths (dashed curves) and constrained paths (solid curves).*

Figure 2 depicts several cumulative occupancy curves $\psi_i$ for a particular example. Suppose that there are $n=100$ types of coupons to collect, and the goal is to collect at least four coupons of as many types as possible (i.e., $I = 3$). Initially, one is given a single coupon of 30 types and pairs of coupons of a further 20 types, corresponding to the initial conditions $\alpha = [0.5\ 0.3\ 0.2]$. In the zero-cost solution (dashed lines), $\psi_3(2.0) \approx 0.71$. Hence, after collecting 200 additional coupons at random ($\beta = 2$), one would expect to have 4 or more coupons for only about 29 types. To compute the likelihood that we have at least 4 coupons for more than 45 types, we take $\xi = 0.55$, which gives a large deviations exponent $J_C \approx 0.18$ and a probability of about $10^{-8}$. The corresponding constrained occupancy curves are depicted by solid lines in Figure 2.

4.4. *Extensions.* There are a number of variations of the basic occupancy problem which can be solved by fairly straightforward generalizations of the results of this paper, but which we will only mention briefly. These include the following problems:

(i) a random number of balls are thrown,
(ii) balls have a probability $p$ of not entering any urn,
(iii) balls enter different subsets of urns with differing probabilities,
(iv) an event of interest may occur at any time in the interval $(0, \beta]$, rather than just at time $\beta$.

Some comments are in order. A particular example of (i) appeared in [13], where the number of balls thrown, $r$, was binomially distributed with parameters $0 < a < 1$ and $n$. An urn model proposed by [16] is of type (ii). Here, the goal is to determine the distribution of the number of targets hit when $r$ shots are fired at $n$ targets, and when the probability of missing the target is $p$. In problem (iii) there are $K > 1$ urn classes with a fraction $\alpha_k, k = 1, 2, \ldots, K$ of urns in each class. Urns enter class $k$ with a fixed probability $p_k$ but then enter any urn within that class with uniform probability. Similar analysis to that for nonempty initial conditions can be applied to this problem. For an example of type (iv), suppose that an infinite sequence of balls is thrown into $n$ initially empty urns, and that one would like to know the probability that the number of urns containing exactly one ball ever exceeds $n/2$. The probability of this occurring *after* $\beta n$ balls have been thrown can be bounded above by the probability that the number of empty and singly occupied urns exceeds $n/2$ when exactly $\beta n$ balls have been thrown. This computation fits into the framework of a partial coupon collectors problem, and the probability can be made negligible from the point of view of large deviations by taking $\beta$ sufficiently large. The remaining possibility, that the event occurs *before* $\beta n$ balls have been thrown, is then a problem of type (iv). The associated calculus of variations problem is to find the lowest cost occupancy curve on $(0, \beta]$ among all curves with $\gamma_1(x) > 0.5$ for some $x \in (0, \beta]$.

## APPENDIX

**Analysis of the calculus of variations problem.** The Appendix is dedicated to proving the calculus of variations results given in Section 2.2. Recall that these results provide explicit representations for the terminal rate function $\mathcal{J}(\omega)$ defined in Corollary 2.3, and for the minimizing occupancy functions $\gamma^*$ satisfying $\mathcal{J}(\omega) = J(\gamma^*)$.

In the first step of the proof, we characterize a set of extremal occupancy paths, that is, paths which satisfy the Euler–Lagrange equations. For all feasible terminal conditions of the form $(1, \omega, \beta)$, Theorem A.1 shows that occupancy paths of the form given in Theorem 2.6 are extremals and that paths of this form can be constructed to meet any feasible constraints of the form $(1, \omega, \beta)$. Likewise, Lemma A.6 and Theorem A.11 show that occupancy paths of the form given in Theorem 2.8 are extremals, and such paths can be constructed to meet all general feasible conditions $(\alpha, \omega, \beta)$. The special form of the extremals is used in Theorem A.5 and Theorem A.13 to show that the extremals have the costs given in Theorem 2.5 and Theorem 2.7, respectively. If $\bar{\gamma}$ is the extremal occupancy path constructed for given constraints $(\alpha, \omega, \beta)$, we thus have the upper bound $\mathcal{J}(\omega) \leq J(\bar{\gamma})$. The assertion in Theorem 2.5 and Theorem 2.7 that the minimum relative entropy is achieved by a unique distribution $\pi^*$ is established in Lemma A.6. The final step needed to prove Theorems 2.5–2.8 is the lower bound $\mathcal{J}(\omega) \geq J(\bar{\gamma})$. This bound is proved in Theorem A.14, using the Euler–Lagrange equations together with properties of the relative entropy.

### A.1. Preliminaries.

A.2. *Proof of Lemma* 2.4. Recall that the lemma states that endpoint constraints $(\alpha, \omega, \beta)$ are feasible if and only if

$$\text{(A.4)} \qquad \sum_{j=0}^{i} \alpha_j \geq \sum_{j=0}^{i} \omega_j, \qquad i = 0, \ldots, I \quad \text{(monotonicity)}$$

and

$$\text{(A.5)} \qquad \sum_{i=0}^{I} i\omega_i + (I+1)\omega_{I+} \leq \sum_{i=0}^{I+1} i\alpha_i + \beta \qquad \text{(conservation)}.$$

PROOF OF LEMMA 2.4. If a valid occupancy curve $\gamma : [0, \beta] \to S_I$ meets the initial and terminal constraints $\alpha$ and $\omega$, then property (b) of Lemma 2.1 implies (A.4) and property (c) implies (A.5) since

$$\sum_{i=0}^{I}\left(\sum_{j=0}^{i}\alpha_j - \sum_{j=0}^{i}\omega_j\right) = \sum_{j=0}^{I}(I+1-j)(\alpha_j - \omega_j)$$

(A.6)

$$= (I+1)(\omega_{I+} - \alpha_{I+1}) + \sum_{j=0}^{I} j(\omega_j - \alpha_j).$$

On the other hand, given the constraints (A.4) and (A.5), one can show that the linear functions

$$\gamma_i(x) = \alpha_i + (\omega_i - \alpha_i)\frac{x}{\beta}, \qquad i = 0, \ldots, I,$$

$$\gamma_{I+}(x) = 1 - \sum_{i=0}^{I}\gamma_i(x) = \alpha_{I+1} + (\omega_{I+} - \alpha_{I+1})\frac{x}{\beta}$$

satisfy the constraints and the conditions of Lemma 2.1. Properties (a) and (b) are immediate, and property (c) will be established by showing that $-\sum_{i=0}^{I}\dot{\psi}_i \leq 1$. Indeed, $-\beta\sum_{i=0}^{I}\dot{\psi}_i$ is equal to the left-hand side of (A.6) and, therefore, $-\sum_{i=0}^{I}\dot{\psi}_i \leq 1$ follows from (A.5). □

A.3. *Euler–Lagrange equations.* Given the numerous descriptions of occupancy processes and rates ($\psi, \gamma, \theta, \dot{\gamma}, \dot{\psi}$, etc.), it is convenient to abuse notation. Thus, for example, we will write both $J(\gamma)$ and $J(\psi)$, with the understanding that the fundamental object of interest is the occupancy function $\gamma$, and that $J(\psi)$ is merely $J(\gamma)$ when $\psi$ is the cumulative occupancy process that corresponds to $\gamma$. Also, since $\dot{\psi}_i = -\theta_i$ for $i = 0, \ldots, I$, we can define the local rate function (which is usually written as a function of $\psi$ and $\dot{\psi}$) as a function of $\psi$ and $\theta$, and represent the overall cost of a cumulative occupancy trajectory $\psi$ as an integral of the form

$$J(\psi) = \int_0^\beta L(\psi(x), \theta(x))\,dx.$$

Because the balls are thrown uniformly and randomly into the urns, the expected rate for balls to enter urns of occupancy $i$ is $\gamma_i = \psi_i - \psi_{i-1}$. As discussed in Section 2, the cost of a deviation of a given path $\psi$ from its expected behavior at a given instant is given by the rate function

$$L(\psi, \theta) = D(\theta\|\gamma)$$

$$= \sum_{i=0}^{I}\theta_i \log\frac{\theta_i}{\gamma_i} + \theta_{I+}\log\frac{\theta_{I+}}{\gamma_{I+}}$$

$$= \sum_{i=0}^{I}\theta_i \log\frac{\theta_i}{\psi_i - \psi_{i-1}} + \left(1 - \sum_{i=0}^{I}\theta_i\right)\log\frac{(1-\sum_{i=0}^{I}\theta_i)}{1-\psi_I}.$$

The rate function is defined to be infinity if the curve is not a cumulative occupancy function.

The calculus of variations problem is to find the path $\psi$ having least cost among all paths satisfying given initial conditions and endpoint constraints, and to find the cost of such a *minimal path*. As illustrated in the examples in Section 4, the results extend to cases where the terminal point (or the initial point) are required to lie in a given constraint set.

DEFINITION A.1. An occupancy path defined on $[0, \beta]$ is said to be an extremal if it satisfies the Euler–Lagrange equations [8],

$$\frac{\partial L}{\partial \psi_i}(\psi(x), \theta(x)) = -\frac{d}{dx}\left\{\frac{\partial L}{\partial \theta_i}(\psi(x), \theta(x))\right\}$$

for all $i \in \{0, \ldots, I\}, x \in (0, \beta)$.

Although the Euler–Lagrange equations are neither necessary nor sufficient conditions for minimality in general, extremals do turn out to be minimal in many cases. In the following sections we will construct a family of extremal paths for the cost function given above, and show that the extremal paths are, in fact, globally minimal.

In the case at hand, the Euler–Lagrange equations are given by

$$(\text{A.7}) \quad -\frac{\theta_i}{\psi_i - \psi_{i-1}} + \frac{\theta_{i+1}}{\psi_{i+1} - \psi_i} = \frac{d}{dx}\left\{-\log\frac{\theta_i}{\psi_i - \psi_{i-1}} + \log\frac{\theta_{I+}}{1 - \psi_I}\right\},$$

for $i = 0, \ldots, I - 1$, and by

$$(\text{A.8}) \quad -\frac{\theta_I}{\psi_I - \psi_{I-1}} + \frac{\theta_{I+}}{1 - \psi_I} = \frac{d}{dx}\left\{-\log\frac{\theta_I}{\psi_I - \psi_{I-1}} + \log\frac{\theta_{I+}}{1 - \psi_I}\right\}.$$

In the case when we have equality in the conservation constraint (2.4), and by taking $I + 1$ to be $I$ if necessary, we have that $\sum_{i=0}^{I} \omega_i = 1$. Such cases of equality are referred to as the polynomial case. When this holds, every valid occupancy path must have $\psi_I(x) = 1$ and, therefore, $\theta_I(x) = \theta_{I+}(x) = 0$. For all such occupancy paths, the rate function $L(\psi, \theta)$ then reduces to

$$(\text{A.9}) \quad L(\psi, \theta) = \sum_{i=0}^{I-1} \theta_i \log \frac{\theta_i}{\psi_i - \psi_{i-1}}.$$

The Euler–Lagrange equations pertinent to the problem of minimizing this restricted set of occupancy paths are just

$$(\text{A.10}) \quad -\frac{\theta_i}{\psi_i - \psi_{i-1}} + \frac{\theta_{i+1}}{\psi_{i+1} - \psi_i} = \frac{d}{dx}\left\{-\log\frac{\theta_i}{\psi_i - \psi_{i-1}}\right\},$$

for $i = 0, \ldots, I - 1$.

**A.4. Characterization of the extremals under empty initial conditions.**
In this section we consider the simplest and most important case, in which the urns are all initially empty, $\alpha_0 = 1$. Recall that to each feasible endpoint constraint, $(1, \omega, \beta)$ correspond to twist parameters $C \geq 0, \rho > 0$ which satisfy the equations

$$\sum_{i=0}^{I} \omega_i + \sum_{i=I+1}^{\infty} C\mathcal{P}_i(\rho\beta) = 1,$$

$$\sum_{i=0}^{I} i\omega_i + \sum_{i=I+1}^{\infty} iC\mathcal{P}_i(\rho\beta) = \beta.$$

THEOREM A.1. *Suppose that $(1, \omega, \beta)$ are feasible terminal constraints, and that $\rho$ and $C$ are the corresponding twist parameters. Then the set of functions $\gamma$ defined by*

(A.8) $$\psi_0(x) = Ce^{-\rho x} + \sum_{k=0}^{I}(\omega_k - C\mathcal{P}_k(\rho\beta))\left(1 - \frac{x}{\beta}\right)^k,$$

(A.9) $$\gamma_i(x) = \frac{x^i}{i!}(-1)^i \psi_0^{(i)}(x), \qquad 0 \leq i \leq I,$$

$$\gamma_{I+}(x) = 1 - \sum_{i=0}^{I} \gamma_i(x),$$

*are extremals on $[0, \beta]$ which satisfy the terminal constraints along with the initial constraint $\psi_0(0) = 1$.*

Recall that in the special case $C = 0$, the first component of the extremal is simply the $I$th order polynomial

(A.10) $$\psi_0(x) = \sum_{k=0}^{I} \omega_k \left(1 - \frac{x}{\beta}\right)^k.$$

We refer to such paths as *polynomial extremals*, and to extremals with $C > 0$ as *exponential extremals*.

The following definition and lemma are useful in the proof of the above theorem.

DEFINITION A.2. *A nonnegative function $\varphi$ is completely monotone on an interval $[a, b]$ if it is infinitely differentiable on $[a, b]$ with*

(A.11) $$(-1)^i \varphi^{(i)}(x) \geq 0 \qquad \text{for all } x \in [a, b] \text{ and } i > 0.$$

This definition is based on the one pertaining to Bernstein's theorem, which characterizes Laplace transforms, see [15]. However, our definition differs in that it considers only a finite interval $[a, b]$.

LEMMA A.2. *The function $\psi_0(x)$ given in (A.8) is completely monotone on $[0, \beta]$. Moreover, the inequality in (A.11) is strict for $x \in [0, \beta)$ and $i = 0, \ldots, I$. In the case $C > 0$, this can be strengthened to all $i = 0, 1, \ldots$.*

PROOF. Clearly, $\psi_0(x)$ is infinitely differentiable. Now suppose first that $C > 0$. The $i$th derivative of $\psi_0$ is

$$(-1)^i \psi_0^{(i)}(x) = \rho^i C e^{-\rho x} + \sum_{k=i}^{I} (\omega_k - C\mathcal{P}_k(\rho\beta)) \beta^{-i} \frac{k!}{(k-i)!} \left(1 - \frac{x}{\beta}\right)^{k-i}$$

for $i = 0, \ldots, I$, and

$$(A.12) \qquad (-1)^i \psi_0^{(i)}(x) = \rho^i C e^{-\rho x}$$

for $i > I$.

It is clear that $a(x) \doteq (-1)^{I+1} \psi_0^{(I+1)}(x)$ is completely monotone on $[0, \infty)$. Moreover, $(-1)^i a^{(i)}(x) > 0$, $x \in [0, \infty)$, $i = 0, 1, \ldots$. We deduce that $\varphi(x) = (-1)^I \psi_0^{(I)}(x)$ must be monotonically strictly decreasing, and that

$$\varphi(\beta) = (-1)^I \psi_0^{(I)}(\beta)$$
$$= \rho^I C e^{-\rho\beta} + \left(\omega_I - C e^{-\rho\beta} \frac{(\rho\beta)^I}{I!}\right) \beta^{-I} I!$$
$$= \left(\frac{\beta^I}{I!}\right)^{-1} \omega_I \geq 0,$$

so that $\varphi(x) > 0$ on $[0, \beta]$. It follows that $\varphi(x)$ is completely monotone on $[0, \beta]$ and that the derivative constraint is strict on $[0, \beta)$. Proceeding inductively to $I - 1$ and beyond, we arrive at the lemma.

In the polynomial case $C = 0$ the argument proceeds similarly on noting that

$$(-1)^I \psi_0^{(I)}(x) = \left(\frac{\beta^I}{I!}\right)^{-1} \omega_I > 0 \qquad \text{for all } x \in [0, \beta]. \qquad \square$$

PROOF OF THEOREM A.1. The terminal constraints can be verified immediately by inspection. The initial constraints follow from the construction of $C$ in (2.6), since

$$\psi_0(0) = C\left(1 - \sum_{i=0}^{I} \mathcal{P}_i(\rho\beta)\right) + \sum_{i=0}^{I} \omega_i = 1.$$

A similar computation also using (2.6) shows that $\psi_0^{(1)}(0) = -1$, a fact that we will need shortly.

To establish that the given functions are valid occupancy curves, it is useful to introduce the infinite sequence of functions $\gamma_i$ and corresponding $\psi_i$ obtained by extending (A.9) to all $i$:

$$(A.13) \qquad \gamma_i(x) = \frac{x^i}{i!}(-1)^i \psi_0^{(i)}(x), \qquad i = 0, 1, \ldots, I, I+1, \ldots.$$

As the sum of an exponential and a polynomial, $\psi_0$ has a Taylor series representation of unlimited radius about any point $x$. Then $\psi_0(0) - \psi_0(x) = \sum_{i=1}^{\infty} \frac{(-x)^i}{i!} \psi_0^{(i)}(x)$, and thus,

$$\sum_{i=0}^{\infty} \gamma_i(x) = \sum_{i=0}^{\infty} \frac{(-x)^i}{i!} \psi_0^{(i)}(x) = \psi_0(0) = 1.$$

Since, by Lemma A.2, the $\gamma_i$ are nonnegative on $[0, \beta]$, we have $\{\gamma_i(x)\} \in S_\infty$ for all $x$ in that interval. It follows from (A.13) that for all $i \geq 0$, the rate of decrease of each cumulative occupancy is

$$(A.14) \begin{aligned} \theta_i(x) &= -\dot{\psi}_i(x) \\ &= -\sum_{k=0}^{i} \dot{\gamma}_k(x) \\ &= \sum_{k=1}^{i} -\frac{x^{k-1}}{(k-1)!}(-1)^k \psi_0^{(k)}(x) + \sum_{k=0}^{i} \frac{x^k}{k!}(-1)^{k+1} \psi_0^{(k+1)}(x) \\ &= \frac{x^i}{i!}(-1)^{i+1} \psi_0^{(i+1)}(x) \geq 0 \end{aligned}$$

for $x \in [0, \beta]$. Forming the Taylor series representation of $\psi_0^{(1)}$ about $x$, it follows from (A.14) that

$$\sum_{i=0}^{\infty} \theta_i(x) = \sum_{i=0}^{\infty} -\dot{\psi}_i(x) = -\psi_0^{(1)}(0) = 1.$$

The infinite sequence of functions can thus be thought of as a valid infinite-dimensional occupancy path on $[0, \beta]$. Conditions (a) and (b) of Lemma 2.1 are immediate from the expressions for $\gamma_i$ and $\dot{\psi}_i$ in terms of $\psi_0^{(i)}$, and condition (c) follows by integrating the inequality

$$\sum_{i=0}^{I} \theta_i(x) \leq 1$$

over an arbitrary subinterval of $[0, \beta]$. Thus, the finite-dimensional occupancy path $\gamma$ is valid.

Finally, we must show that the given curves solve the Euler–Lagrange differential equations in $(0, \beta)$. We begin with the exponential case $C > 0$. The $I+$ terms satisfy the simple expressions

$$\theta_{I+}(x) = 1 + \sum_{i=0}^{I} \dot{\psi}_i(x) = \sum_{i=I+1}^{\infty} -\dot{\psi}_i(x) = \rho C \sum_{i=I+1}^{\infty} \mathcal{P}_i(\rho x),$$

$$1 - \psi_I(x) = \sum_{i=I+1}^{\infty} \gamma_i(x) = C \sum_{i=I+1}^{\infty} \mathcal{P}_i(\rho x),$$

where the first display uses $\theta = -\dot{\psi}$ and equations (A.12) and (A.14), and the second uses (A.12) and (A.13). Then the ratio $\theta_{I+}(x)/(1 - \psi_I(x)) = \rho$ is constant, and the corresponding terms drop out of the right-hand side of (A.7) and (A.8).

From the expressions for $\gamma_i(x)$ and $\dot{\psi}_i(x)$ given in (A.13) and (A.14), it follows that

$$(A.15) \qquad \frac{\theta_i(x)}{\gamma_i(x)} = -\frac{\psi_0^{(i+1)}(x)}{\psi_0^{(i)}(x)}, \qquad x \in (0, \beta), i \geq 0.$$

Then

$$(A.16) \quad \frac{d}{dx}\left\{-\log \frac{\theta_i}{\gamma_i}\right\} = -\frac{\psi_0^{(i+2)}}{\psi_0^{(i+1)}} + \frac{\psi_0^{(i+1)}}{\psi_0^{(i)}} = \frac{\theta_{i+1}}{\gamma_{i+1}} - \frac{\theta_i}{\gamma_i}, \qquad x \in (0, \beta),$$

verifying (A.7). Likewise, to verify (A.8), we apply (A.16) with $i = I$, using the substitution $-\psi_0^{(I+2)}/\psi_0^{(I+1)} = \rho$ obtained from (A.12).

In the polynomial case $C = 0$, note that (A.15) holds for $i = 0, \ldots, I$, so that (A.16) implies (A.10). □

A.5. *Interpretation of the twist parameter $\rho$.* For an exponential extremal the twist parameter $\rho$ may be interpreted as follows. The expected rate for balls to enter urns with more than $I$ balls is equal to the proportion of such urns, namely, $1 - \psi_I$. The twist parameter $\rho$ is then a multiplicative factor applied to $1 - \psi_I$ to give the actual rate at which balls enter these urns. Thus, if $\rho > 1$, balls unusually pile into high-occupancy urns, while if $\rho < 1$, they instead concentrate on the low-occupancy urns. The occupancy distribution of the high-occupancy urns remains Poisson but with a modified parameter.

The infinite sequence of occupancy functions introduced in the proof just given is a useful construct. Operations which are technically difficult in infinite dimensions may nevertheless be carried out *formally* on the infinite sequence of functions, giving insight in to the solution for the finite-dimensional system.

The next two lemmas compute the cost of extremal curves in a general form which will also apply to the nonempty case of the next section. The conditions of the lemmas are satisfied by the exponential and polynomial extremals of Theorem A.1, as can be readily verified.

LEMMA A.3. *Suppose that $\psi_0$ is completely monotone on $[0, \beta]$ with*

$$(-1)^i \psi_0^{(i)}(x) = C\rho^i e^{-\rho x} \qquad \text{for } i > I,$$

*for some $I \geq 0$, $C > 0$, and $\rho > 0$. Further suppose that $\{\gamma_0, \gamma_1, \ldots\}$ are an infinite sequence of nonnegative functions on $[0, \beta]$ satisfying*

$$\sum_{i=0}^{\infty} \gamma_i(x) = 1, \qquad \sum_{i=0}^{\infty} i\gamma_i(x) \leq B_0,$$

$$\sum_{i=0}^{\infty} -\dot{\psi}_i(x) = 1, \qquad -\dot{\psi}_i(x) = \gamma_i(x) \frac{-\psi_0^{(i+1)}(x)}{\psi_0^{(i)}(x)},$$

*for all $x \in [0, \beta]$ and some constant $B_0 < \infty$. Let $\gamma$ denote the vector of functions $\{\gamma_0, \ldots, \gamma_I, \gamma_{I+}\}$, where $\gamma_{I+} = 1 - \sum_{i=0}^{I} \gamma_i$. Then the cost $J(\gamma)$ is given by*

$$(A.17) \qquad J(\gamma) = \beta + \sum_{i=0}^{\infty} [\gamma_i(\beta) \log|\psi_0^{(i)}(\beta)| - \gamma_i(0) \log|\psi_0^{(i)}(0)|].$$

PROOF. For all $i > I$ and $x \in [0, \beta]$, we have $\psi_0^{(i+1)}/\psi_0^{(i)} = -\rho$. Therefore, using $\theta = -\dot{\psi}$ with (A.12) and (A.13), it follows that

$$\theta_{I+} = 1 - \sum_{i=0}^{I} -\dot{\psi}_i = \sum_{i=I+1}^{\infty} -\dot{\psi}_i = \rho \sum_{i=I+1}^{\infty} \gamma_i = \rho(1 - \psi_I).$$

Then for each $x \in [0, \beta]$, the cost function $L$ can be interpreted as an infinite-dimensional cost function $L_\infty$. Indeed, since $-\dot{\psi}_i/\gamma_i = -\psi_0^{(i+1)}/\psi_0^{(i)} = \rho$ for $i > I$ and $(1 + \sum_{i=0}^{I} \dot{\psi}_i)/(1 - \psi_I) = \rho$,

$$L(\psi, \dot{\psi}) = \sum_{i=0}^{I} (-\dot{\psi}_i) \log \frac{(-\dot{\psi}_i)}{\psi_i - \psi_{i-1}}$$

$$+ (1 + \dot{\psi}_0 + \cdots + \dot{\psi}_I) \log \frac{(1 + \dot{\psi}_0 + \cdots + \dot{\psi}_I)}{1 - \psi_I}$$

$$= \sum_{i=0}^{\infty} (-\dot{\psi}_i) \log \frac{(-\dot{\psi}_i)}{\psi_i - \psi_{i-1}} \equiv L_\infty(\psi, \dot{\psi}).$$

Note that since $L_\infty(\psi, \dot\psi)$ can be interpreted as a relative entropy, we always have $L_\infty(\psi, \dot\psi) \geq 0$. The total cost $J(\gamma)$ may be computed by integrating $L_\infty$. Note that $|\sum_i \dot\psi_i(x)| = 1, x \in [0, \beta]$ and that given $\varepsilon > 0, \log -\psi_0^{(i+1)}(x)/\psi_0^{(i)}(x)$ are uniformly bounded for $x \in [\varepsilon, \beta - \varepsilon]$. By the monotone convergence theorem

$$J(\gamma) = \lim_{\varepsilon \downarrow 0} \int_\varepsilon^{\beta-\varepsilon} \sum_{i=0}^\infty -\dot\psi_i(x) \log \frac{-\psi_0^{(i+1)}(x)}{\psi_0^{(i)}(x)} \, dx.$$

Using our convention that $\psi_{-1}(x) = 0$, the integral on the right may be written as

$$\int_\varepsilon^{\beta-\varepsilon} \sum_{i=0}^\infty -\dot\psi_i(x) \log|\psi_0^{(i+1)}(x)| \, dx + \int_\varepsilon^{\beta-\varepsilon} \sum_{i=0}^\infty \dot\psi_i(x) \log|\psi_0^{(i)}(x)| \, dx$$

(A.18)
$$= \sum_{i=0}^\infty \int_\varepsilon^{\beta-\varepsilon} (\dot\psi_i(x) - \dot\psi_{i-1}(x)) \log|\psi_0^{(i)}(x)|.$$

The above expression is valid as long as the left and right series in the first line converge. But this follows from the finite mean condition on $\{\gamma_i\}(x)$ since, for $i > I$,

$$-\dot\psi_i(x) \log|\psi_0^{(i)}(x)| = \rho \gamma_i(x)[\log C + i \log \rho - \rho x]$$

and similarly for $\dot\psi_i(x) \log |\psi_0^{(i+1)}(x)|$. Applying integration by parts and using (A.15) for each term of (A.18), the integral is

$$\sum_{i=0}^\infty [\gamma_i \log|\psi_0^{(i)}|]_\varepsilon^{\beta-\varepsilon} + \int_\varepsilon^{\beta-\varepsilon} \sum_{i=0}^\infty -\dot\psi_i \, dx.$$

The lemma follows on taking limits as $\varepsilon \downarrow 0$. □

A similar result holds in the case of the polynomial extremal.

LEMMA A.4.  *Suppose that $\psi_0$ is a degree $I$ polynomial which is completely monotone on $[0, \beta]$. Let $\{\gamma_0, \ldots, \gamma_I\}$ be nonnegative functions on $[0, \beta]$ satisfying*

$$\sum_{i=0}^I \gamma_i(x) = 1, \qquad \sum_{i=0}^{I-1} -\dot\psi_i(x) = 1,$$

$$-\dot\psi_i(x) = \gamma_i(x) \frac{-\psi_0^{(i+1)}(x)}{\psi_0^{(i)}(x)}, \qquad i = \{0, \ldots, I-1\}$$

*for all $x \in [0, \beta]$. Then the cost $J(\gamma)$ is given by*

(A.19)  $$J(\gamma) = \beta + \sum_{i=0}^I [\gamma_i(\beta) \log|\psi_0^{(i)}(\beta)| - \gamma_i(0) \log|\psi_0^{(i)}(0)|].$$

PROOF. The cost obtained by integrating the reduced cost function (A.9) from 0 to $\beta$, using the same substitutions and integration by parts as in the proof of Theorem A.3. □

A.6. *Characterization of the extremal cost.* In the case of empty initial conditions, recall that the $\gamma_i$ defined in the proof of Theorem A.1 satisfy $\gamma_i(x) = x^i |\psi_0^{(i)}(x)|/i!$. Using this expression to substitute for $\psi_0^{(i)}$ in (A.17) and (A.19), we find that the cost $J(\gamma)$ is simply the relative entropy $D(\gamma(\beta)\|\mathcal{P}(\beta))$ between the $\gamma_i(\beta)$ and the Poisson, zero cost distribution. It turns out that the given $\gamma_i(\beta)$ minimize the relative entropy, among all distributions for which the first $I+1$ elements are determined by $\omega$, and which have mean $\beta$. Denoting the set of all such distributions by $F(1, \omega, \beta)$, we may prove Theorem 2.5, which we restate here.

THEOREM A.5. *Suppose that $\gamma$ is an extremal occupancy path constructed according to Theorem A.1 to meet feasible terminal constraints $(1, \omega, \beta)$. Then*

$$J(\gamma) = \min_{\pi \in F(1,\omega,\beta)} D(\pi \| \mathcal{P}(\beta)).$$

PROOF. We first solve the minimization problem, and then relate the problem to $J(\gamma)$. The result is trivially true for polynomial extremals, since then $F(1, \omega, \beta)$ has only one element. In the case of an exponential extremal, the given minimization problem can be solved using Lagrange multipliers, which turn out to be simple functions of the twist parameters $\rho$ and $C > 0$ associated with the given endpoint constraints. We consider the set of nonnegative sequences $\{\pi_j\}_{j=0}^\infty \in \mathbb{R}_+^\infty$ satisfying $\pi_i = \omega_i$ for $i = 0, \ldots, I$, and define the *Lagrangian* $\mathcal{L}$ on this set to be

$$\mathcal{L}(\pi; y, z) = \sum_{i=0}^\infty \pi_i \log \frac{\pi_i}{\mathcal{P}_i(\beta)} + y\left(1 - \sum_{i=0}^\infty \pi_i\right) + z\left(\beta - \sum_{i=0}^\infty i\pi_i\right),$$

where $z \doteq \log \rho$ and $y \doteq \log C + (1-\rho)\beta + 1$.

Define $\pi_i^* \doteq \gamma_i(\beta)$, where the latter are determined as in the proof of Theorem A.1 with the given $C, \rho$. For any $i > I$, the definitions of $\pi_i^*, x$ and $y$ and the strict convexity of $x \log x$ imply that $\pi_i^*$ is the unique global minimizer of $x \to x \log x - x[\log \mathcal{P}_i(\beta) + y + iz]$. Therefore, for any $i > I$ and $\pi_i \in [0, \infty)$,

$$\pi_i \log \frac{\pi_i}{\mathcal{P}_i(\beta)} - y\pi_i - zi\pi_i \geq \pi_i^* \log \frac{\pi_i^*}{\mathcal{P}_i(\beta)} - y\pi_i^* - zi\pi_i^*.$$

Following standard Lagrangian arguments, we thus have

$$\inf_{\pi \in F(1,\omega,\beta)} D(\pi \| \mathcal{P}(\beta)) = \inf_{\tilde{\pi} \in F(1,\omega,\beta)} \mathcal{L}(\tilde{\pi}; y, z)$$
$$\geq \inf_{\pi \in \mathbb{R}_+^\infty} \mathcal{L}(\pi; y, z)$$
$$= D(\pi^* \| \mathcal{P}(\beta)).$$

Since $\pi^* \in F(1,\omega,\beta)$, it follows that $\pi^*$, the terminal distribution, is the minimizer of the relative entropy. The uniqueness of $\pi^*$ follows from the strict convexity of the relative entropy with respect to its first argument.

Substituting the particular form of $\pi_i^*$ into (A.17), we have

$$J(\gamma) = \sum_{i=0}^{\infty} \gamma_i(\beta) \log \frac{\gamma_i(\beta)}{\beta^i e^{-\beta}/i!} = D(\pi^* \| \mathcal{P}(\beta)),$$

where

$$\gamma_i(\beta) = \begin{cases} \omega_i, & 0 \leq i \leq I, \\ C\mathcal{P}_i(\rho\beta), & i > I. \end{cases} \qquad \square$$

The cost for an exponential extremal can be written explicitly in terms of $\rho$, $\omega$ and $\beta$ as

$$J(\gamma) = \sum_{i=0}^{I} \omega_i \log \frac{\omega_i}{\mathcal{P}_i(\beta)} + \left(1 - \sum_{i=0}^{I} \omega_i\right)(\log C + (1-\rho)\beta)$$
$$+ \left(\beta - \sum_{i=0}^{I} i\omega_i\right) \log \rho,$$

where $C$ is defined by (2.6). In the polynomial case, the cost is simply given by the first term in the above expression.

**A.7. Characterization of the extremals under general initial conditions.**
We now generalize the solution to the Euler–Lagrange equations for the case when the urns are not all initially empty, but in which they may contain up to $K$ balls. Recall that the fraction of urns having $k$ balls initially is denoted $\alpha_k$, and the set of urn occupancies $k$ with $\alpha_k > 0$ is denoted $\mathcal{K}$. Without loss of generality, we may assume that some urns are initially empty (hence, $0 \in \mathcal{K}$) and we may take $K = \max \mathcal{K}$.

In the case $K \geq 1$, we regard the urns as belonging to $|\mathcal{K}|$ classes, according to their initial occupancies. We denote the final number of additional balls per urn entering the $k$th class by $\beta_k$, in which case the total number of additional balls per urn is $\beta = \sum_{k \in \mathcal{K}} \alpha_k \beta_k$. It turns out in the solution of the Euler–Lagrange equations that the fraction of balls entering the $k$th class in any time period is $\alpha_k \beta_k / \beta$. After rescaling time by the factor $\beta_k/\beta$, the evolution of *additional balls* entering each class of urns becomes an occupancy problem with initially empty conditions and terminal time $\beta_k$.

We can thus define $|\mathcal{K}|$ occupancy curves $\gamma_{(k)} = \{\gamma_{k,i}\}_i$, where the function $\gamma_{k,i} : [0, \beta_k] \to [0, 1]$ represents the fraction of class $k$ urns which contain $i$ additional balls (thus, $k+i$ total balls) after $x$ balls per urn have been given to this class. The overall extremal occupancy curves for the general initial conditions are then given by

$$\gamma_i(x) = \sum_{k \in \mathcal{K}} \alpha_k \gamma_{k,i-k}(x\beta_k/\beta).$$

We will see that the occupancy curve $\gamma_{(k)}$ for the $k$th subproblem is an extremal of the form given by Theorem A.1, with $I - k$ terminal constraints. Given the appropriate subproblem terminal conditions $(1, \omega_{(k)}, \beta_k)$, the results of the previous section determine the twist parameters $C_k$ and $\rho_k$ and corresponding extremals.

At this point we come to the main obstacle in determining the extremals for nonempty initial conditions. This is to show that there *exist* subproblem terminal conditions $(1, \omega_{(k)}, \beta_k)$ which yield the extremal solution in the overall problem. To show that such conditions and the corresponding extremals exist, we will first give the form of the final cost function (which can be obtained by formal arguments). We then show that the cost function is the solution of a minimization problem, that the problem has a unique minimizing argument, that it has corresponding Lagrange multipliers, and that the extremal curves can be constructed using the Lagrange multipliers.

The large deviations exponent for the case of nonempty initial conditions $\alpha$ turns out to be

(A.20) $$\min_{\pi \in S_\infty^{|\mathcal{K}|}} \sum_{k \in \mathcal{K}} \alpha_k D(\pi_{(k)} \| \mathcal{P}(\beta)),$$

subject to the conservation constraint

(A.21) $$\sum_{k \in \mathcal{K}} \alpha_k \sum_{j=0}^{\infty} j \pi_{k,j} = \beta$$

and the terminal conditions

(A.22) $$\omega_i = \sum_{k \leq i, k \in \mathcal{K}} \alpha_k \pi_{k,i-k} \quad \text{for all } 0 \leq i \leq I.$$

Recall that $S_\infty^{|\mathcal{K}|}$ is the set of all $|\mathcal{K}|$-tuples of distributions $\{\pi_{(k)}\}_{k \in \mathcal{K}}$, where each $\pi_{(k)}$ is a distribution on the nonnegative integers. It is straightforward to show that this problem is feasible whenever $(\alpha, \omega, \beta)$ are feasible constraints, as we discuss in the proof of Lemma A.6. As in the empty case, each polynomial problem can be formulated in a standard way so that $\omega_{I+} = 0$, in which case equality holds in the monotonicity condition $\sum_{j=0}^{I} \alpha_j = \sum_{j=0}^{I} \omega_j = 1$ as well as in the conservation condition (2.4). Unlike in the empty case, the minimization problem does not become trivial in the polynomial case since degrees of freedom remain in allocating balls among the $|\mathcal{K}|$ classes. In the polynomial case, the constraints imply that

$$\pi_{k,j} = 0, \qquad k+j > I \tag{A.23}$$

so that the minimization problem is finite dimensional. The polynomial problem can be stated equivalently as minimization of (A.20) subject to (A.23) and the endpoint constraints (A.22), in which case the conservation constraint (A.21) need not be included explicitly.

As written, the vector of distributions $\pi = (\pi_{(0)}, \ldots, \pi_{(k)}, \ldots, \pi_{(K)})$ is such that $k$ only ranges over indices with $\alpha_k > 0$. Including all indices $0 \leq k \leq K$ yields an equivalent problem in which the minimizing solution is not unique, since the $\pi_{(k)}$ with $\alpha_k = 0$ make no contribution to the objective function or the constraints. In the form given, however, the solution can be shown to be unique.

LEMMA A.6. *Suppose that $(\alpha, \omega, \beta)$ are feasible constraints. Then there is a unique vector of distributions $\pi^* \in S_\infty^{|\mathcal{K}|}$ which minimize (A.20) subject to constraints (A.21) and (A.22).*

PROOF. Recall that $F(\alpha, \omega, \beta)$ denotes the set of all distributions $\pi \in S_\infty^{|\mathcal{K}|}$ satisfying the constraints (A.21) and (A.22). The feasibility of $(\alpha, \omega, \beta)$ implies that $F$ has at least one element. Moreover, it has at least one element with finite support and, thus, finite cost. A sketch argument for this point is as follows: First of all, any set of exponential constraints with terminal condition $\omega \in S_I$ can be reformulated as a set of polynomial constraints with terminal condition $\hat{\omega} \in S_{\hat{I}}$, where $\hat{I} > I$. Any feasible point for constraints $(\alpha, \hat{\omega}, \beta)$ would also be feasible for the original constraints, and would have finite support. It remains to show that there is at least one feasible point for each set of feasible polynomial constraints. Such a feasible point may be constructed by an ordered filling construction. First, we assign $\pi_{0,0} = \omega_0/\alpha_0$. The first monotonicity condition ($i = 0$) in (2.3) ensures that this is possible with less than or equal to unit mass. Next, some of the remaining mass from the distribution $\pi_{(0)}$ is applied to $\pi_{0,1}$, and mass from $\pi_{(1)}$ is applied to $\pi_{1,0}$ as necessary until the constraint $\alpha_0 \pi_{0,1} + \alpha_1 \pi_{1,0} = \omega_1$ is satisfied. That $\pi_{(0)}$ and $\pi_{(1)}$ have sufficient mass to do so follows from the ($i = 1$) condition in (2.3). This process continues until all constraints $\omega_i$ have been satisfied. The equality in the conservation condition (2.4) implies that the final step uses up all of the probability mass in the distributions $\{\pi_{(k)}\}$.

It will be useful to consider $|\mathcal{K}|$-tuples of distributions in $S_\infty$ under the topology of weak convergence, with the distance between distributions given by the Prohorov metric [5]. For two distributions $P$, $Q$ in $S_\infty$, this distance is simply

$$\sum_{i:P_i>Q_i} P_i - Q_i = \sum_{i:Q_i>P_i} Q_i - P_i,$$

and the metric extends to $|\mathcal{K}|$-tuples by treating them as elements of a product space. For each $A < \infty$, define the set

$$H_A = \left\{\pi \in S^{|\mathcal{K}|} : \sum_k \alpha_k D(\pi_{(k)} \| \mathcal{P}(\beta)) \leq A\right\}.$$

The level sets of the relative entropy are compact under the above topology ([11], Lemma 1.4.3), from which it follows that $H_A$ is also compact. For later reference, we note that a finite sum of relative entropy functions also inherits from the relative entropy the properties of lower semi-continuity and of strict convexity with respect to its first argument. As the next step in proving the lemma, we wish to show for some finite $A$ that

$$Q_A(\alpha, \omega, \beta) \doteq F(\alpha, \omega, \beta) \cap H_A$$

is compact and nonempty. This will enable us to find the minimum as a limit of a sequence of distributions in $Q_A$.

Since there are solutions with finite cost, it is automatic that $Q_A(\alpha, \omega, \beta)$ is nonempty for large enough $A$. Since the set $Q_A$ is bounded, it is enough to verify for each convergent sequence of distributions $\pi^{(n)} \in Q_A$ that the limiting distribution $\overline{\pi}$ lies in $F(\alpha, \beta, \omega)$. (That it lies in $H_A$ is immediate since this set is compact.) The only difficulty lies in showing that the mean of the $k$th subdistribution $\overline{\pi}_{(k)}$ is equal to the limit of the means of the $\pi_{(k)}^{(n)}$. This will be the case if we can show for any sequence in $Q_A$ and for an arbitrary $k$ that the sequence $\pi_{(k)}^{(n)}$ is uniformly integrable. In the present context, uniform integrability means that for each $\varepsilon > 0$, there is $m < \infty$ such that

(A.24) $$\sum_{j \geq m} j \pi_{k,j}^{(n)} \leq \varepsilon$$

for all $n$ and $k$.

Since $e^y - 1 = \sup_{x>0}[xy - x\log x + x - 1]$, we have the inequality

$$xy \leq (x \log x - x + 1) + (e^y - 1),$$

where by convention $0 \log 0 = 0$. Observe that this inequality is valid for all $y$ and all $x \geq 0$, and that $x \log x - x + 1$ is nonnegative. We will also make use of the fact that the Poisson distribution has exponential moments: for any $\delta > 0$,

$$\sum_{j=0}^\infty e^{1/\delta j} \mathcal{P}_j(\beta) = \sum_{j=0}^\infty \frac{e^{1/\delta j} e^{-\beta} \beta^j}{j!} = e^{\beta e^{1/\delta} - \beta} < \infty.$$

Taking $y \doteq j/\delta$ and $x \doteq \pi_{k,j}^{(n)}/\mathcal{P}_j(\beta)$, we have the estimate

$$\sum_{j \geq m} j\pi_{k,j}^{(n)} = \sum_{j \geq m} \delta \frac{j}{\delta} \frac{\pi_{k,j}^{(n)}}{\mathcal{P}_j(\beta)} \mathcal{P}_j(\beta)$$

$$\leq \sum_{j \geq 0} \delta \left[ \pi_{k,j}^{(n)} \log\left(\frac{\pi_{k,j}^{(n)}}{\mathcal{P}_j(\beta)}\right) - \pi_{k,j}^{(n)} + \mathcal{P}_j(\beta) \right] + \delta \sum_{j \geq m} (e^{1/\delta j} - 1)\mathcal{P}_j(\beta)$$

$$= \delta D(\pi_k^{(n)} \| \mathcal{P}(\beta)) + \delta \sum_{j \geq m} (e^{1/\delta j} - 1)\mathcal{P}_j(\beta)$$

$$\leq \delta A + \delta \sum_{j \geq m} (e^{1/\delta j} - 1)\mathcal{P}_j(\beta),$$

where the second equality uses the fact that both $\pi_{(k)}^{(n)}$ and $\mathcal{P}(\beta)$ are probability distributions. Thus, (A.24) follows by first picking $\delta > 0$ sufficiently small and then $m < \infty$ sufficiently large. [Note also for use in Corollary A.7 that the analogous result holds if $\beta$ is replaced by a sequence $\beta^{(n)} \to \beta \in (0, \infty)$].

We now apply the uniform integrability to analyze the limits of the means. It follows from [5], Theorem 5.4, that

$$\beta_k^{(n)} \to \overline{\beta}_k.$$

Finally,

$$\sum_k \alpha_k \overline{\beta}_k = \lim_n \sum_k \alpha_k \beta_k^{(n)} = \lim_n \beta = \beta,$$

so that $\overline{\pi} \in F(\beta, \alpha, \omega)$ and $Q_A$ is compact.

We have shown that $Q_A$ is compact under the Prohorov metric and that it is nonempty. Now let

$$G \doteq \inf_{\pi \in Q_A} \sum_{k \in \mathcal{K}} \alpha_k D(\pi_{(k)} \| \mathcal{P}(\beta)) \geq 0.$$

Choose a sequence $\pi^{(n)} \in Q_A$ such that $\sum_{k=0}^{K} \alpha_k D(\pi_{(k)}^{(n)} \| \mathcal{P}(\beta)) \leq G + 1/n$. Since $Q_A$ is compact this has a convergent subsequence, and to simplify the notation we index this subsequence by $n$. Let $\pi^* \in Q_A$ be the limit point. Since relative entropy $D(\pi \| \eta)$ is lower semi-continuous in $\pi$ ([11], Lemma 1.4.3(b)), it follows that

$$\sum_{k \in \mathcal{K}} \alpha_k D(\pi_{(k)}^* \| \mathcal{P}(\beta)) \leq \liminf_n \sum_{k \in \mathcal{K}} \alpha_k D(\pi_{(k)}^{(n)} \| \mathcal{P}(\beta)) = G,$$

where in fact we have equality by definition of $G$. The uniqueness of this minimizing vector of distributions follows from the strict convexity of the objective function. $\square$

Before continuing with the question of the existence of the extremals, we pause to obtain a corollary which will be useful in establishing the strong minimum in Theorem A.14. Define $\mathcal{J}(\alpha,\omega,\beta)$ to be the minimum cost in the problem (A.20) if the constraints are feasible, and $\infty$ otherwise.

COROLLARY A.7. *Suppose that $(\alpha^{(n)},\omega^{(n)},\beta^{(n)})$ is a sequence of feasible optimization problems with costs $\mathcal{J}^{(n)} \doteq \mathcal{J}(\alpha^{(n)},\omega^{(n)},\beta^{(n)})$ such that*

$$(\alpha^{(n)},\omega^{(n)},\beta^{(n)}) \to (\alpha,\omega,\beta)$$

*componentwise and such that $(\alpha,\omega,\beta)$ is feasible with $0<\beta<\infty$. Then*

$$\liminf_n \mathcal{J}^{(n)} \geq \mathcal{J}(\alpha,\omega,\beta).$$

PROOF. For any $A > \liminf_n \mathcal{J}^{(n)} \doteq G$, we may choose a subsequence of problems such that

$$\mathcal{J}^{(n)} \leq G + 1/n < A,$$

and for each problem in this subsequence, we define $\pi^{(n)} \in Q_A(\alpha^{(n)},\beta^{(n)},\omega^{(n)})$ to be the minimizing solution $\pi^*$ given in the above lemma. By compactness, there is a further subsequence such that $\pi^{(n)} \to \overline{\pi} \in H_A$. As in the proof of Lemma A.6, if $\beta < \infty$, then

$$\beta_k^{(n)} \doteq \sum_{j=0}^{\infty} j\pi_{k,j}^{(n)} \to \beta_k \doteq \sum_{j=0}^{\infty} j\overline{\pi}_{k,j}, \qquad k=0,1,\ldots,K.$$

If $\alpha_k^{(n)} \to 0$ and $\beta_k^{(n)}\alpha_k^{(n)} \not\to 0$, then the term $\alpha_k^{(n)} D(\pi_{(k)}^{(n)}\|\mathcal{P}(\beta^{(n)}))$ would go to infinity; to see this, note that the infimum of $D(\pi\|\mathcal{P}(\beta))$ over distributions $\pi$ with mean $\lambda$ is $\beta - \lambda + \lambda\log(\lambda/\beta)$. Since our sequence has bounded cost, we must have $\alpha_k^{(n)}\beta_k^{(n)} \to \alpha_k\beta_k$ even if $\alpha_k = 0$.

Thus,

$$\sum_{k\in\mathcal{K}} \beta_k\alpha_k = \lim_n \sum_k \beta_k^{(n)}\alpha_k^{(n)} = \lim_n \beta^{(n)} = \beta,$$

and $\overline{\pi}$ satisfies the constraints $(\alpha,\omega,\beta)$. By joint lower semi-continuity of the relative entropy in both arguments ([11], Lemma 1.4.3(b)) we have

$$\liminf_n \sum_{k=0}^{K} \alpha_k^{(n)} D(\pi_{(k)}^{(n)}\|\mathcal{P}(\beta^{(n)})) \geq \liminf_n \sum_{k\,:\,\alpha_k>0} \alpha_k^{(n)} D(\pi_{(k)}^{(n)}\|\mathcal{P}(\beta^{(n)}))$$
$$\geq \sum_{k\,:\,\alpha_k>0} \alpha_k D(\overline{\pi}_{(k)}\|\mathcal{P}(\beta))$$
$$\geq \mathcal{J}(\alpha,\omega,\beta)$$

as desired. □

We now turn to the existence of the Lagrange multipliers corresponding to the minimization problem (A.20). For a given vector $\omega$, it will be helpful to define the set of integers $\mathcal{I}(\omega)$, where $i \in \mathcal{I}(\omega)$ if $0 \leq i \leq I$ and $\omega_i > 0$, or if $i > I$ and $\omega_{I+} > 0$. This is the set of terminal urn occupancy levels which the constraints do not force to be empty. As we will show, for irreducible constraints, the optimal $\{\pi_{k,j}\}$ with $k+j \in \mathcal{I}(\omega)$ are always strictly positive. This result does not hold directly for reducible constraints, but as we discussed earlier, any problem with reducible constraints can be replaced by a finite number of subproblems with irreducible constraints.

LEMMA A.8. *Suppose that $(\alpha, \omega, \beta)$ are irreducible feasible constraints, and let $\pi^* \in S_\infty^{|\mathcal{K}|}$ be the minimizer of (A.20) subject to constraints (A.21) and (A.22). Then $\pi_{k,j}^* > 0$ for all $k+j \in \mathcal{I}(\omega)$ and $\pi_{k,j}^* = 0$ if $k+j \notin \mathcal{I}(\omega)$.*

PROOF. For $k+j \notin \mathcal{I}(\omega)$, the constraints force the $\pi_{k,j}$ to be zero for all feasible points of the problem, and hence, for the minimizer in particular. The main point is to show that it is feasible for any other element to be positive, and the result will then follow from the infinite derivative of the objective function near the boundary.

Let $\pi \in S_\infty^{|\mathcal{K}|}$ be any feasible solution meeting the given constraints, and suppose that $\pi_{m,n} = 0$ for some particular $(m,n)$ with $m+n \in \mathcal{I}(\omega)$. Let $\tilde{\pi}$ be an arbitrary set of probability distributions in $S_\infty^{|\mathcal{K}|}$ subject to the restrictions that $\tilde{\pi}_{m,n} > 0$, that $\tilde{\pi}_{k,j} = 0$ for all $k+j \notin \mathcal{I}(\omega)$, and that $\tilde{\pi}$ has finite support. If we define

$$\tilde{\omega}_i \doteq \sum_{k=0}^{i} \alpha_k \tilde{\pi}_{k,i-k}, \qquad \tilde{\beta} \doteq \sum_{k \in \mathcal{K}} \alpha_k \sum_{j=0}^{\infty} j \tilde{\pi}_{k,j},$$

then $\tilde{\pi}$ is a feasible solution for the constraints $(\alpha, \tilde{\omega}, \tilde{\beta})$. Next, for any $\eta > 0$, we may define

$$\hat{\beta} \doteq (\beta - \eta\tilde{\beta})/(1-\eta), \qquad \hat{\omega} \doteq (\omega - \eta\tilde{\omega})/(1-\eta).$$

By construction, the constraints $(\alpha, \hat{\omega}, \hat{\beta})$ are feasible for sufficiently small $\eta$. In the exponential case, this is immediately true because $(\alpha, \omega, \beta)$ satisfies all constraints in (2.3) and (2.4) with strict inequality. In the polynomial case, the $(\alpha, \omega, \beta)$ satisfies the $I$th monotonicity constraint of (2.3) and the conservation constraint (2.4) with equality. However, it may easily be verified that these two constraints also hold with equality for $(\alpha, \tilde{\omega}, \tilde{\beta})$ and, hence, for $(\alpha, \hat{\omega}, \hat{\beta})$.

We may now let $\hat{\pi}$ be a feasible point for $(\alpha, \hat{\omega}, \hat{\beta})$, and finally form $\bar{\pi} = \eta \tilde{\pi} + (1-\eta)\hat{\pi}$. By construction, $\bar{\pi}$ is a feasible solution corresponding to the original constraints $(\alpha, \omega, \beta)$, and, in addition, $\bar{\pi}_{m,n} > 0$. As discussed in the proof of Lemma A.6, we may take $\hat{\pi}$ to have finite support so that $\bar{\pi}$ also has finite support.

We will show that $\pi$ is not a minimizer by proving that a sufficiently small perturbation toward $\bar{\pi}$ reduces the objective function. Because the constraints are linear, the points $\pi^\varepsilon = \varepsilon \bar{\pi} + (1-\varepsilon)\pi$ are feasible solutions to the original constraints for all $0 \leq \varepsilon \leq 1$. Let $f(\varepsilon)$ denote the value of the objective function evaluated at $\pi^\varepsilon$. The derivative of this function is

$$\begin{aligned}
f'(\varepsilon) &= \sum_{k \in \mathcal{K}} \alpha_k \sum_{j=0}^{\infty} \left( \log \frac{\pi^\varepsilon_{k,j}}{\mathcal{P}_j(\beta)} + 1 \right) (\bar{\pi}_{k,j} - \pi_{k,j}) \\
&= \sum_{k \in \mathcal{K}} \alpha_k \sum_{j=0}^{\infty} \log \frac{\pi^\varepsilon_{k,j}}{\mathcal{P}_j(\beta)} (\bar{\pi}_{k,j} - \pi_{k,j}) \\
&= \alpha_m \bar{\pi}_{m,n} \log \frac{\pi^\varepsilon_{m,n}}{\mathcal{P}_n(\beta)} + \sum_{k \in \mathcal{K}} \alpha_k \sum_{(k,j) \neq (m,n)} \bar{\pi}_{k,j} \log \frac{\pi^\varepsilon_{k,j}}{\mathcal{P}_j(\beta)} \\
&\quad - \sum_{k \in \mathcal{K}} \alpha_k \sum_j \pi_{k,j} \log \frac{\pi^\varepsilon_{k,j}}{\mathcal{P}_j(\beta)}.
\end{aligned}$$

In the limit as $\varepsilon \to 0$, the third term in the last display tends to

$$-\sum_k \alpha_k D(\pi_{(k)} \| \mathcal{P}(\beta)) \leq 0.$$

The second term is bounded above by the expression $-\sum_k \alpha_k \sum_j \bar{\pi}_{k,j} \log \mathcal{P}_j(\beta)$, which is finite because $\bar{\pi}$ has finite support. The first term in the display tends to $-\infty$, which establishes that $f(\delta) < f(0)$ for some sufficiently small $\delta > 0$. We have shown that, for any $m+n \in \mathcal{I}(\omega)$, $\pi$ cannot minimize (A.20) if $\pi_{m,n} = 0$. As $\pi_{m,n}$ is arbitrary within $m+n \in \mathcal{I}(\omega)$, it follows that $\pi^*_{m,n} > 0$ whenever $m+n \in \mathcal{I}(\omega)$. $\square$

We now establish the existence of Lagrange multipliers and the form of the optimal solution for the polynomial case.

THEOREM A.9. *Suppose that $(\alpha, \omega, \beta)$ are irreducible feasible constraints yielding equality in* (2.4), *and let $\pi^*$ be the corresponding unique minimizer from Lemma* A.6. *Then there exist positive constants $\{D_k\}_{k \in \mathcal{K}}$ and $\{W_i\}_{i \in \mathcal{I}}$ such that the minimizer takes the form*

$$\pi^*_{k,j} = \begin{cases} D_k \mathcal{P}_j(\beta) W_{k+j}, & k \in \mathcal{K},\, j+k \in \mathcal{I}(\omega), \\ 0, & \text{otherwise.} \end{cases}$$

PROOF. The strict equality in the constraints requires every feasible point to be supported on the finite set satisfying $k \in \mathcal{K}$ and $k + j \in \mathcal{I}$ and, hence, the minimization problem (A.20) may be considered to be a finite-dimensional problem over this set. By Lemma A.8 the minimizer $\pi^*$ is strictly positive. Since the objective function is continuously differentiable in a neighborhood of $\pi^*$, and the constraints are linear, Lagrange multipliers are guaranteed to exist for this problem (see, e.g., [4], Proposition 1.33). Specifically, there are constants $z_k$ and $w_i$ such that the Lagrangian

$$\mathcal{L}(\pi) = \sum_{k \in \mathcal{K}} \alpha_k \sum_{k+j \in \mathcal{I}} \pi_{k,j} \log \frac{\pi_{k,j}}{\mathcal{P}_j(\beta)} + \sum_{k \in \mathcal{K}} z_k \alpha_k \left(1 - \sum_{k+l \in \mathcal{I}} \pi_{k,l}\right) + \sum_{i \in \mathcal{I}} w_i \left(\omega_i - \sum_{k=0}^{i} \alpha_k \pi_{k,i-k}\right)$$

is also minimized at $\pi^*$, with all partial derivatives of $\mathcal{L}$ being zero at the optimal point. Taking partial derivatives and rearranging, the optimality condition yields

$$\pi^*_{k,j} = \mathcal{P}_j(\beta) e^{z_k - 1 + w_{k+j}}.$$

The result follows on defining $D_k = e^{z_k - 1}$ and $W_i = e^{w_i}$. □

We now give the corresponding theorem for the exponential case.

THEOREM A.10. *Suppose that $(\alpha, \omega, \beta)$ are irreducible feasible constraints yielding strict inequality in* (2.4), *and let $\pi^*$ be the corresponding unique minimizer from Lemma* A.6. *Then there exist positive constants $\rho$, $\{C_k\}_{k \in \mathcal{K}}$, and $\{W_i\}_{i \in \mathcal{I}, i \leq I}$ such that the minimizer takes the form*

$$\pi^*_{k,j} = \begin{cases} C_k \mathcal{P}_j(\rho\beta) W_{k+j}, & k \in \mathcal{K},\, k+j \leq I,\, k+j \in \mathcal{I}, \\ C_k \mathcal{P}_j(\rho\beta), & k \in \mathcal{K},\, k+j > I, \\ 0, & \text{otherwise.} \end{cases}$$

PROOF. To avoid difficulties with infinite-dimensional Lagrangians, consider a sequence of truncated problems indexed by a sequence of integers $M > I$. For each such $M$ define

$$\beta^{(M)} \doteq \sum_{k \in \mathcal{K}} \alpha_k \sum_{l \leq M} l \pi^*_{k,l}, \qquad \eta_k^{(M)} \doteq \sum_{l \leq M} \pi^*_{k,l}$$

and consider the problem of minimizing

$$\sum_{k \in \mathcal{K}} \alpha_k \sum_{j \leq M, k+j \in \mathcal{I}} \pi_{k,j} \log \frac{\pi_{k,j}}{\mathcal{P}_j(\beta)},$$

subject to the constraints

$$\sum_{k \in \mathcal{K}} \alpha_k \sum_{l \leq M, k+l \in \mathcal{I}} l \pi_{k,l} = \beta^{(M)},$$

$$\sum_{l \leq M, k+l \in \mathcal{I}} \pi_{k,l} = \eta_k^{(M)}, \qquad k \in \mathcal{K},$$

$$\sum_{k \leq i, k \in \mathcal{K}} \alpha_k \pi_{k,i-k} = \omega_i, \qquad i \in \mathcal{I}.$$

By construction, the minimizer of this problem is obtained simply by truncating $\pi^*$. As in the previous lemma, the strict positivity of the minimizer, together with the linearity of the constraints, guarantees the existence of Lagrange multipliers via [4], Proposition 1.33. The Lagrangian for the $M$th problem is

$$\mathcal{L}^{(M)} \doteq \sum_{k \in \mathcal{K}} \alpha_k \sum_{j \leq M, k+j \in \mathcal{I}} \pi_{k,j} \log \frac{\pi_{k,j}}{\mathcal{P}_j(\beta)} + y^{(M)} \left( \beta^{(M)} - \sum_{k \in \mathcal{K}} \alpha_k \sum_{l \leq M, k+l \in \mathcal{I}} l \pi_{k,l} \right)$$
$$+ \sum_{k \in \mathcal{K}}^K z_k^{(M)} \alpha_k \left( \eta_k^{(M)} - \sum_{l \leq M, k+l \in \mathcal{I}} \pi_{k,l} \right)$$
$$+ \sum_{i \in \mathcal{I}} w_i^{(M)} \left( \omega_i - \sum_{k \leq i, k \in \mathcal{K}} \alpha_k \pi_{k,i-k} \right).$$

Since the derivative of $\mathcal{L}$ with respect to $\pi_{k,j}$ must be zero at the minimizer, it follows that

$$\pi_{k,j}^* = \mathcal{P}_j(\beta) e^{jy^{(M)} - 1 + z_k^{(M)} + w_{k+j}^{(M)}},$$

for all $j < M$ and $k + j \in \mathcal{I}$, where for convenience we have defined $w_i = 0$ for $i > I$. For $k + j > I$, note that

$$\frac{\pi_{k,j+1}^*}{\pi_{k,j}^*} = \frac{\mathcal{P}_{j+1}(\beta)}{\mathcal{P}_j(\beta)} e^{y^{(M)}}$$

so that $y^{(M)}$ is independent of $M$. Since $w_{k+j}^{(M)} = 0$ if $k + j > I$, it follows that $z_k^{(M)}$ and $w_i^{(M)}$ are also independent of $M$. Then the above expression for $\pi_{k,j}^*$ holds for all $k + j \in \mathcal{I}$, with fixed Lagrange multipliers $y$, $\{z_k\}$ and $\{w_i\}$. The form expressed in the theorem is based on the substitutions $\rho = e^y$, $C_k = e^{(\rho-1)\beta - 1 + z_k}$ and $W_i = e^{w_i}$. □

Now that we have the general form of the minimizing endpoint values $\pi_{k,j}^*$, we are ready to characterize the extremal curves. We denote the mean of the $k$th distribution $\pi_{(k)}^*$ by $\beta_k \doteq \sum_j j \pi_{k,j}^*$. In the exponential case, we define $\rho_k = \rho\beta/\beta_k$ so that the distribution $\pi_{(k)}$ meets $I - k$ terminal constraints, has mean $\beta_k$ and a Poisson $\rho_k\beta_k$ tail. Denoting the terminal constraints by $\omega_{k,j} \doteq \pi_{k,j}^*$ and $\omega_{(k)} = (\omega_{k,0}, \ldots, \omega_{k,I-k})$, we see that $\pi_{(k)}^*$ is the minimizing distribution arising in Theorem A.1 for the terminal constraints $(1, \omega_{(k)}, \beta_k)$, with associated twist parameters $\rho_k$ and $C_k$. In the polynomial case, the terminal constraints $(1, \omega_{(k)}, \beta_k)$ mean that each subproblem is polynomial, with $C_k = 0$. In either case, Theorem A.1 gives the form of the extremal occupancy curves $\gamma_{k,j}(x)$ for each of these $k$ subproblems. We now show that the subproblem curves sum to form the general extremals.

THEOREM A.11. *Suppose that $\{\pi_{k,j}^*\}$ are the minimizing distributions from Lemma A.6 for feasible constraints $(\alpha, \omega, \beta)$. Denote the means of the minimizing arguments by $\beta_k = \sum_j j\pi_{k,j}^*$, and let $\gamma_{k,j}(x)$ be the extremal curves from Theorem A.1 for the subproblems $(1, \omega_{(k)}, \beta_k)$. Then the curves*

$$\gamma_i(x) = \sum_{k=0}^{i} \alpha_k \gamma_{k,i-k}(x\beta_k/\beta), \qquad i = 0, \ldots, I,$$

$$\gamma_{I+}(x) = 1 - \sum_{i=0}^{I} \gamma_i(x)$$

*are extremals which satisfy the constraints $(\alpha, \omega, \beta)$.*

PROOF. We assume without loss of generality that the constraints are irreducible. Otherwise, each irreducible subproblem may be treated separately.

We extend the definition of $\gamma_i(x)$ for $i > I$ by using the extended definition of $\gamma_{k,j}(x)$ used in the proof of Theorem A.1 [see (A.13)]. Also, for convenience, define $\overline{\gamma}_{k,j}(x) \doteq \gamma_{k,j}(x\beta_k/\beta)$ and likewise, define $\overline{\psi}_{k,j}$. The $\psi_i$ inherit from the $\gamma_i$ the relation

$$\psi_i(x) = \sum_{k=0}^{i} \alpha_k \overline{\psi}_{k,i-k}(x).$$

To see that the $\gamma_i$ are valid occupancy curves, note that

$$\sum_{i=0}^{\infty} \gamma_i(x) = \sum_{i=0}^{\infty} \sum_{k=0}^{i} \alpha_k \overline{\gamma}_{k,i-k}(x) = \sum_{k=0}^{K} \alpha_k \sum_{j=0}^{\infty} \overline{\gamma}_{k,j}(x) = 1.$$

Also, the $-\dot{\psi}_i(x)$ are all nonnegative on $[0, \beta]$, and

$$\sum_{i=0}^{\infty} -\dot{\psi}_i(x) = \sum_{k=0}^{K} \alpha_k \sum_{j=0}^{\infty} -\dot{\overline{\psi}}_{k,j}(x) = \sum_{k=0}^{K} \alpha_k(\beta_k/\beta) = 1,$$

where the last equality follows from the fact that the $\pi^*_{k,j}$ satisfy the conservation constraint (A.21). Hence, curves $\gamma_i(x), i \leq I$ and $\gamma_{I+}(x)$ satisfy the conditions of Lemma 2.1.

It is clear that the $\gamma_i$ satisfy the desired initial conditions since $\gamma_{k,0}(0) = 1$ and $\gamma_{k,j}(0) = 0$ for $j > 0$. The terminal conditions are guaranteed by the fact that the values $\gamma_{k,j}(\beta_k) = \pi^*_{k,j}$ satisfy the constraints (A.22).

In order to establish that the curves satisfy the Euler–Lagrange equations, we first show that the rescaled zero-occupancy curve $\overline{\psi}_{k,0}$ for the $k$th subproblem is proportional to the $k$th derivative of the overall zero-occupancy curve $\psi_0$.

First, take the exponential case, and recall that $\rho_k \beta_k = \rho \beta$. The $k$th zero-occupancy curve, after rescaling, is

$$\overline{\psi}_{k,0}(x) = C_k e^{-\rho_k(x\beta_k/\beta)} + \sum_{i=0}^{I-k}(\omega_{k,i} - C_k \mathcal{P}_i(\rho_k \beta_k))\left(1 - \frac{x\beta_k/\beta}{\beta_k}\right)^i$$

$$= C_k e^{-\rho x} + \sum_{i=0}^{I-k}(\omega_{k,i} - C_k \mathcal{P}_i(\rho\beta))\left(1 - \frac{x}{\beta}\right)^i$$

$$= C_k\left[e^{-\rho x} + \sum_{i=0}^{I-k}(W_{k+i} - 1)\mathcal{P}_i(\rho\beta)\left(1 - \frac{x}{\beta}\right)^i\right].$$

Here we have used that $\omega_{k,i} = \pi^*_{k,i} = C_k \mathcal{P}_i(\beta) W_{k+i}$ when $i \leq I - k$. The $k$th derivative of $\psi_0 = \alpha_0 \overline{\psi}_{0,0}$ is

(A.25)
$$(-1)^k \psi_0^{(k)}(x)$$
$$= \alpha_0 C_0 \left[\rho^k e^{-\rho x} + \sum_{i=k}^{I}(W_i - 1)\mathcal{P}_i(\rho\beta)\frac{i!}{(i-k)!\beta^k}\left(1 - \frac{x}{\beta}\right)^{i-k}\right]$$
$$= \alpha_0 C_0 \rho^k \left[e^{-\rho x} + \sum_{j=0}^{I-k}(W_{k+j} - 1)\mathcal{P}_j(\rho\beta)\left(1 - \frac{x}{\beta}\right)^j\right]$$
$$= \frac{\alpha_0 C_0 \rho^k}{C_k}\overline{\psi}_{k,0}(x).$$

Using Theorem A.1 to express $\gamma_{k,j}$ in terms of $\psi_{k,0}$, we have

$$\gamma_i(x) = \sum_{k=0}^{i} \alpha_k \frac{(x\beta_k/\beta)^{i-k}}{(i-k)!}(-1)^{i-k}\psi_{k,0}^{(i-k)}\left(\frac{x\beta_k}{\beta}\right)$$

$$= \sum_{k=0}^{i} \alpha_k \frac{x^{i-k}}{(i-k)!}(-1)^{i-k}\frac{d^{i-k}}{dx^{i-k}}\overline{\psi}_{k,0}(x)$$

$$= \sum_{k=0}^{i} \frac{\alpha_k C_k}{\alpha_0 C_0 \rho^k} \frac{x^{i-k}}{(i-k)!} (-1)^i \psi_0^{(i)}(x).$$

Similarly, we obtain

$$-\dot{\psi}_i(x) = -\sum_{k=0}^{i} \alpha_k \frac{\beta_k}{\beta} \dot{\psi}_{k,i-k}\left(\frac{x\beta_k}{\beta}\right)$$

$$= \sum_{k=0}^{i} \alpha_k \frac{x^{i-k}}{(i-k)!} (-1)^{i-k} \frac{d^{i-k+1}}{dx^{i-k+1}} \overline{\psi}_{k,0}(x)$$

$$= \sum_{k=0}^{i} \frac{\alpha_k C_k}{\alpha_0 C_0 \rho^k} \frac{x^{i-k}}{(i-k)!} (-1)^{i+1} \psi_0^{(i+1)}(x),$$

so that the extremals satisfy the simple relation

(A.26) $$\frac{\theta_i(x)}{\gamma_i(x)} = -\frac{\psi_0^{(i+1)}(x)}{\psi_0^{(i)}(x)},$$

which also arose in the case of empty initial conditions. Because the polynomial portion of $\psi_0$ has degree $I$, the ratio $\theta_i(x)/\gamma_i(x)$ is given by the constant $\rho$ for $i > I$, which establishes (as in the proof of Lemma A.3) that

(A.27) $$\frac{\theta_{I+}(x)}{1 - \psi_I(x)} = \rho.$$

As shown in the proof of Theorem A.1, the relations (A.26) and (A.27) are sufficient to show that the $\gamma_i$ satisfy the Euler–Lagrange equations.

In the polynomial case, similar computations show that

(A.28) $$(-1)^k \psi_0^{(k)}(x) = \frac{\alpha_0 D_0}{D_k} \overline{\psi}_{k,0}(x)$$

and that

(A.29) $$\theta_i(x) = \sum_{k=0}^{i} \frac{\alpha_k D_k}{\alpha_0 D_0} \frac{x^{i-k}}{(i-k)!} (-1)^{i+1} \psi_0^{(i+1)}(x) = \gamma_i(x) \frac{-\psi_0^{(i+1)}(x)}{\psi_0^{(i)}(x)}$$

for $i = 0, \ldots, I$, establishing the restricted set of equations (A.10). □

In order to demonstrate a strong minimum we will need the following.

COROLLARY A.12. *For irreducible constraints, the occupancy functions defined in Theorem A.11 satisfy the integrated version of the Euler–Lagrange equations: given $x_1 \in [0, \beta)$, there are constants $C_i, i = 0, \ldots, I-1$ (and also $i = I$ in the exponential case) depending on $x_1$ only such that*

$$\int_{x_1}^{x'} \frac{\partial L}{\partial \psi_i}(\psi(x), \theta(x))\, dx + \frac{\partial L}{\partial \theta_i}(\psi(x'), \theta(x')) = C_i$$

*for all $x' \in (x_1, \beta)$.*

PROOF. We may take $\alpha_0 > 0$ without loss of generality. By Theorem A.11, the extremals satisfy the usual form of the Euler–Lagrange equations. The above indefinite integral is $\log(-\psi_0^{(i+1)}(x)/\psi_0^{(i)}(x))$, which is finite at both endpoints, because $(-1)^i \psi_0^{(i)}(x) > 0$ for $x \in [0, \beta)$, $i \leq I$ (and $i = I + 1$ in the exponential case), as shown in Lemma A.2. The partial derivative with respect to $\theta_i$ also exists for the same reason. □

THEOREM A.13. *Suppose that $\gamma$ is the extremal defined by Theorem A.11 for the feasible constraints $(\alpha, \omega, \beta)$. Then the cost $J(\gamma)$ is the solution to the minimization problem (A.20) subject to constraints (A.21) and (A.22).*

PROOF. Consider first the exponential case. The infinite sequence of functions $\gamma_i$ defined in the proof of Theorem A.11 are shown in the proof to satisfy the conditions of Lemma A.3. Hence, the cost is

$$J(\gamma) = \beta + \sum_{i=0}^{\infty}\left[\sum_{k=0}^{i} \alpha_k \overline{\gamma}_{k,i-k}(\beta) \log|\psi_0^{(i)}(\beta)|\right] - \sum_{k=0}^{K} \alpha_k \log|\psi_0^{(k)}(0)|$$

$$= \sum_{k=0}^{K} \alpha_k \sum_{j=0}^{\infty} \overline{\gamma}_{k,j}(\beta) \log\left|\frac{\psi_0^{(k+j)}(\beta)}{\psi_0^{(k)}(0) e^{-\beta}}\right|.$$

The fact that $\psi_{k,0}(0) = 1$ together with (A.25) implies that $\psi_0^{(k)}(x) = \psi_0^{(k)}(0) \times \overline{\psi}_{k,0}(x)$, so that

$$\frac{\psi_0^{(k+j)}(x)}{\psi_0^{(k)}(0) e^{-\beta}} = \left(\frac{\beta_k}{\beta}\right)^j \psi_{k,0}^{(j)}\left(\frac{x\beta_k}{\beta}\right) e^{\beta} = \frac{\overline{\gamma}_{k,j}(x)}{(-x)^j e^{-\beta}/j!}.$$

Then

$$J(\gamma) = \sum_{k=0}^{K} \alpha_k D(\overline{\gamma}_{k,\cdot}(\beta) \| \mathcal{P}(\beta)).$$

By construction, the endpoints $\overline{\gamma}_{k,j}(\beta)$ coincide with the optimal arguments $\pi_{k,j}^*$ of the minimization problem.

The proof of the polynomial case is almost identical, except that we use Lemma A.4 and (A.28). □

**A.8. Extremal curves have globally minimal cost.** In this section we prove the following theorem:

THEOREM A.14 (Strong minimum). *Given feasible constraints $(\alpha, \omega, \beta)$, let $\gamma$ be the corresponding extremal occupancy path defined in Theorem A.11, and let $\tilde{\gamma}$ be any other occupancy path satisfying the same constraints. Then $J(\gamma) \leq J(\tilde{\gamma})$.*

We first introduce some notation. Let $\mathcal{O}$ denote the set of valid occupancy functions, that is, vector functions $\gamma$ such that the cumulative occupancy functions $\psi$ satisfy the conditions of Lemma 2.1, and let $\mathcal{O}(\alpha, \omega, \beta) \doteq \{\gamma \in \mathcal{O} : \gamma(0) = \alpha, \gamma(\beta) = \omega\}$ be the subset of valid occupancy functions satisfying feasible constraints $(\alpha, \omega, \beta)$. The proof of the following lemma is a straightforward consequence of the convexity of the map $(\theta, \gamma) \to D(\theta \| \gamma)$ and, hence, omitted.

LEMMA A.15.

(a) $\mathcal{O}(\alpha, \omega, \beta)$ *is a convex set.*
(b) $J(\gamma)$ *restricted to $\mathcal{O}(\alpha, \omega, \beta)$ is a convex function.*

For a given pair of occupancy functions $\bar{\gamma}, \tilde{\gamma} \in \mathcal{O}(\alpha, \omega, \beta)$, we denote $\gamma^\varepsilon \doteq (1-\varepsilon)\bar{\gamma} + \varepsilon\tilde{\gamma}$ and define the function $G : [0,1] \to \mathbb{R}_+ \cup \{\infty\}$ by

$$G[\varepsilon] \doteq J(\gamma^\varepsilon) = \int_0^\beta D(\theta^\varepsilon(x) \| \gamma^\varepsilon(x)) \, dx.$$

To simplify the notation, we do not explicitly indicate the dependence of $G$ on $\bar{\gamma}$ and $\tilde{\gamma}$. Lemma A.15 ensures that $G$ is a well-defined, convex function for any $\bar{\gamma}, \tilde{\gamma} \in \mathcal{O}(\alpha, \omega, \beta)$.

For the remainder of this section, we once again restrict attention to irreducible constraints, without loss of generality. The following lemma establishes the minimality of the extremals in an important special case. The proof uses the fact that the extremal curve and its derivative can be bounded away from zero.

LEMMA A.16. *Suppose that $(\alpha, \omega, \beta)$ are strictly positive, irreducible feasible constraints, where the upper indices $K$ and $I$ of $\alpha$ and $\omega$ satisfy $K = I+1$ in the exponential case, or $K = I$ in the polynomial case. Suppose that $\gamma$ is the extremal curve for these constraints, defined in Theorem A.11, and that $\tilde{\gamma}$ is any competing occupancy function satisfying the same constraints. Then*

$$J(\gamma) \leq J(\tilde{\gamma}).$$

PROOF. We construct the family of paths $\gamma^\varepsilon = (1-\varepsilon)\gamma + \varepsilon\tilde{\gamma}$, with $G[\varepsilon] = J(\gamma^\varepsilon)$. It follows from convexity that $G$ is left and right differentiable wherever it is finite. We will show that $G'_+[0] = 0$, where $G'_+[\varepsilon]$ denotes the right derivative of $G$. The convexity of $G$ then implies the desired result.

It will be convenient to work with the cumulative occupancy functions $\psi^\varepsilon$, and as in Section A.3, we will mix notation by writing, for example, $J(\gamma) = \int_0^\beta L(\psi, \theta) \, dx$.

After defining
$$\eta(x) \doteq \tilde{\psi}(x) - \psi(x),$$
we may write

(A.30) $$G[\varepsilon] = \int_0^\beta L(\psi + \varepsilon\eta, \theta - \varepsilon\dot{\eta}) \, dx \doteq \int_0^\beta g(x, \varepsilon) \, dx.$$

We can assume without loss that $G(1) = J(\tilde{\gamma}) < \infty$, since there is nothing to prove otherwise. We wish to show that differentiation under the integral sign with respect to $\varepsilon$ is valid in a neighborhood of 0. The validity of this operation will follow from [23], Corollary 39.2, if we can provide a constant bound on the partial derivative of $g$ with respect to $\varepsilon$ for almost every $x \in [0, \beta]$.

To construct this bound, we will first establish that the components $\gamma_i(x)$ and derivatives $\theta_i(x) = -\dot{\psi}_i(x)$ are uniformly bounded away from zero. For specificity, we first assume that $(\alpha, \omega, \beta)$ are exponential constraints. Note that in the empty case studied in Section A.4, $\gamma_0(x)$ decreases monotonically to $\gamma_0(\beta) = \omega_0$ and $\theta_0(x)$ decreases monotonically to $\theta_0(\beta)$. Inspection of the formula for $\psi_0$ in Theorem A.1 reveals that $\theta_0(\beta) = \omega_1/\beta$ if $I > 0$ and $\theta_0(\beta) = C\mathcal{P}_1(\rho\beta)/\beta$ otherwise. Hence, for any exponential problem with empty initial conditions and positive terminal conditions $\omega_i$, $\gamma_0$ and $\theta_0$ are uniformly bounded away from zero. For the constraints $(\alpha, \omega, \beta)$ under consideration, each of the associated subproblems is an exponential empty problem with positive terminal conditions, and so $\gamma_{k,0}(x)$ and $\theta_{k,0}(x)$ are uniformly bounded from zero for each subproblem. These, in turn, lower bound the overall $\gamma_i$ and $\theta_i$, since, for example,
$$\gamma_i(x) = \sum_{k=0}^i \alpha_k \gamma_{k,i-k}(x) \geq \alpha_i \gamma_{i,0}(x).$$

The catch-all term $\gamma_{I+}(x)$ is monotonically increasing, and hence, satisfies $\gamma_{I+}(x) \geq \alpha_{I+1} > 0$. The identity $\theta_{I+}(x) = \rho\gamma_{I+}(x)$, established in the proof of Theorem A.11, ensures a similar bound on $\theta_{I+}(x)$.

Note that $\gamma^\varepsilon \geq (1-\varepsilon)\gamma$ and $\theta^\varepsilon \geq (1-\varepsilon)\theta$, so that for each $0 \leq \varepsilon < 1$, these functions are also uniformly bounded from zero on $[0, \beta]$. Indeed, given an arbitrary $\varepsilon_0 < 1$, there are positive lower bounds on $\gamma^\varepsilon$ and $\theta^\varepsilon$ which hold uniformly for $x \in [0, \beta]$ and $\varepsilon \in [0, \varepsilon_0]$. To be precise, these inequalities and bounds hold everywhere except possibly on a set of measure zero where $\tilde{\theta}$ may not exist.

The partial derivative of the integrand of (A.30) is

$$\frac{\partial}{\partial \varepsilon} g(x, \varepsilon) = \sum_{i=0}^{I} \left.\frac{\partial L}{\partial \psi_i}\right|_{\gamma^\varepsilon} (x) \eta_i(x) - \sum_{i=0}^{I} \left.\frac{\partial L}{\partial \theta_i}\right|_{\gamma^\varepsilon} (x) \dot{\eta}_i(x).$$

The partial derivatives of $L$ are (in mixed notation)

$$\frac{\partial L}{\partial \psi_i}(x) = -\frac{\theta_i(x)}{\gamma_i(x)} + \frac{\theta_{i+1}(x)}{\gamma_{i+1}(x)},$$

$$\frac{\partial L}{\partial \theta_i}(x) = \log \frac{\theta_i(x)}{\gamma_i(x)} - \log \frac{\theta_{I+}(x)}{1 - \psi_I(x)}.$$

The uniform lower bounds on $\gamma^\varepsilon$ and $\theta^\varepsilon$, together with upper bounds $\gamma^\varepsilon \leq 1$, $\theta^\varepsilon \leq 1$ and the boundedness of $\eta$ and $\dot{\eta}$, combine to establish that there is a finite $B$ such that

$$\left|\frac{\partial}{\partial \varepsilon} g(x, \varepsilon)\right| < B \qquad \text{for } \varepsilon \in [0, \varepsilon_0], \text{ a.e. } x \in [0, \beta].$$

Similar arguments may be used to establish this bound in the polynomial case. In that case, the terms $\psi_I = 1$ and $\theta_I = 0$ do not play a role, and the expression for the partial derivative of $L$ with respect to $\theta_i$ simplifies.

We are now free to differentiate under the integral sign so that

$$G'[\varepsilon] = \int_0^\beta \frac{\partial}{\partial \varepsilon} g(x, \varepsilon) \, dx$$

for all $\varepsilon \in [0, \varepsilon_0)$. We note that

$$\int_0^x \left.\frac{\partial L}{\partial \psi_i}\right|_\gamma dx' = \int_0^x \left\{-\frac{\theta_i}{\gamma_i} + \frac{\theta_{i+1}}{\gamma_{i+1}}\right\}\bigg|_\gamma dx'$$

is absolutely continuous as a function of $x$, since it is the difference of two continuous monotone functions with bounded derivatives. Applying integration by parts for absolutely continuous functions ([23], 36.1, page 209) to the first term, we obtain

$$G'_+[0] = \int_0^\beta \left\{\frac{\partial L}{\partial \psi}(x) \cdot \eta(x) - \frac{\partial L}{\partial \theta}(x) \cdot \dot{\eta}(x)\right\} dx$$

$$= \int_0^\beta \dot{\eta}(x) \cdot \left(-\int_0^x \frac{\partial L}{\partial \psi}(x') \, dx' - \frac{\partial L}{\partial \theta}(x)\right) dx$$

$$= -\int_0^\beta \sum_{i=0}^I C_i \dot{\eta}_i(x) \, dx$$

$$= \sum_{i=0}^I C_i(\eta_i(0) - \eta_i(\beta)) = 0.$$

The constants $C_i$ appearing in the third equality come from the integrated version of the Euler–Lagrange equations (Corollary A.12), while the last equality is due to the fact that $\gamma$ and $\tilde{\gamma}$ have the same beginning and end points. $\square$

We now extend the lemma to show that the extremal is a global minimizer even when some of the initial and terminal points are zero. Without loss of generality we suppose that $\alpha_0 > 0$ and note that the extremals defined in Theorem A.1 and in Theorem A.11 are strictly positive except possibly at the initial and final times $t = 0$ and $t = \beta$.

Let $\gamma$ be this extremal and suppose that $\tilde{\gamma}$ is an alternate occupancy function. It may be supposed that $\tilde{\gamma}$ also lies on the boundary only at $t = 0$ or $t = \beta$, since if there is a $\tilde{\gamma}$ with lower cost than $\gamma$, the convexity of $J$ implies that there is another occupancy path of the form $\lambda\tilde{\gamma} + (1-\lambda)\gamma$ which also has lower cost than $\gamma$, and which avoids the boundary everywhere that $\gamma$ does.

Given $x_1^{(n)} \downarrow 0$ and $x_2^{(n)} \uparrow \beta$, let $\gamma^{(n)}$ be the extremal curve with initial point $\tilde{\gamma}(x_1^{(n)})$ and terminal point $\tilde{\gamma}(x_2^{(n)})$. By Lemma A.16,

$$J(\gamma^{(n)}) \leq \int_{x_1^{(n)}}^{x_2^{(n)}} D(\tilde{\theta}(x) \| \tilde{\gamma}(x))\, dx.$$

It then follows that

$$\begin{aligned} J(\gamma) &= \mathcal{J}(\alpha, \omega, \beta) \\ &\leq \liminf_n \mathcal{J}(\tilde{\gamma}(x_1^{(n)}), \tilde{\gamma}(x_2^{(n)}), x_2^{(n)} - x_1^{(n)}) \\ &= \liminf_n J(\gamma^{(n)}) \\ &\leq \lim_n \int_{x_1^{(n)}}^{x_2^{(n)}} D(\tilde{\theta}(x) \| \tilde{\gamma}(x))\, dx \\ &= J(\tilde{\gamma}) \end{aligned}$$

with the first two equalities following from Theorem A.13 and the definition of $\mathcal{J}(\alpha, \omega, \beta)$, the first inequality following from Corollary A.7 and the last equality from the monotone convergence theorem.

**Acknowledgments.** We would like to thank an associate editor and two referees for suggestions that significantly improved this paper.

LARGE DEVIATION FOR OCCUPANCY PROBLEMS 57[23] McShane, E. J. (1947). *Integration*. Princeton Univ. Press.
[24] Ramakrishna, M. V. and Mukhopadhyay, P. (1988). Analysis of bounded disorder file organisation. In *Proc. of the 7th ACM Symposium on Principles of DataBase Systems, Austin, Texas* 117–125.
[25] Sagan, H. (1969). *Introduction to the Calculus of Variations*. Dover, New York. MR1210325
[26] Shwartz, A. and Weiss, A. (1995). *Large Deviations for Performance Analysis*. Chapman and Hall, New York. MR1335456
[27] Shwartz, A. and Weiss, A. (2002). Large deviations with diminishing rates. Unpublished manuscript.
[28] Weiss, A. (1994). Large deviations for the occupancy problem. Private communication.
P. Dupuis
Lefschetz Center for Dynamical Studies
Brown University
Providence, Rhode Island 02912
USA
e-mail: dupuis@dam.brown.edu
url: cm.bell-labs.com/who/nuzman

C. Nuzman
P. Whiting
Bell Labs
600 Mountain Avenue
Murray Hill, New Jersey 07974
USA
e-mail: cnuzman@lucent.com
e-mail: pawhiting@lucent.com
url: cm.bell-labs.com/who/pawhiting
e-mail: **??**